# A ROBUST METHOD FOR CLUSTER ANALYSIS

By María Teresa Gallegos and Gunter Ritter

*Universität Passau*

Let there be given a contaminated list of $n$ $\mathbb{R}^d$-valued observations coming from $g$ different, normally distributed populations with a common covariance matrix. We compute the ML-estimator with respect to a certain statistical model with $n-r$ outliers for the parameters of the $g$ populations; it detects outliers and simultaneously partitions their complement into $g$ clusters. It turns out that the estimator unites both the minimum-covariance-determinant rejection method and the well-known pooled determinant criterion of cluster analysis. We also propose an efficient algorithm for approximating this estimator and study its breakdown points for mean values and pooled SSP matrix.

**1. Introduction.** The aim of cluster analysis is the partitioning of a data set into $g$ disjoint subsets or clusters with common characteristics. Besides *heuristics*, there are important approaches to this problem based on *statistical models*, in particular, approaches by the ML and Bayesian paradigms. The latter offer several advantages. They allow one to *compute* the cluster criteria to be optimized and they yield *algorithms* that effectively, and sometimes efficiently, reduce them; see Schroeder (1976). Finally, a model serves as a guide for the user in which cases to apply the method.

This paper deals with statistical cluster analysis in the potential presence of *contaminations*. Statistical methods postulate that the data come from different statistical populations. After clustering, the elements of the clusters may be used in order to estimate the parameters of the underlying statistical laws. Since almost all real data contain outliers, for the method to be useful in practice one will have to allow that part of the data are contaminations or spurious elements. Accommodating or discarding them in a previous step is necessary for robustly estimating these parameters.

There are a great number of statistical techniques for the clustering problem $g \geq 2$ in the *absence* of outliers. One distinguishes between *mixture* and *classification* models; for an overview see Hartigan (1975) and the recent review paper Fraley and Raftery (2002). Two (nowadays classical) statistical partitioning methods are the trace and determinant criteria of cluster analysis [Friedman and Rubin (1967) and Scott and Symons (1971)]. In both criteria, the pooled within-groups sum of squares and products (SSP) matrix **W** of the clustering, see (2), plays a central role. These criteria postulate as estimators those partitions which minimize the trace and the determinant of **W**, respectively. Both methods are not only heuristically motivated: the resulting partitions are maximum likelihood estimators of well-defined statistical models. Therefore, both methods perform well whenever the data set is a realization of random variables obeying the underlying statistical laws. The probabilistic model for which the trace criterion is optimal assumes that all populations are normally distributed with *unknown* mean vectors and the same *spherical* covariance matrix of unknown size. The determinant criterion retains the assumption on equality of the covariance matrices, but is less restrictive in dropping that on sphericity. As a consequence, the partition optimal for the determinant criterion is invariant not only w.r.t. location but also w.r.t. scale.

In the *presence* of outliers and in the case $g = 1$, the problem reduces to outlier detection or robust parameter estimation and a great number of methods are available; for a good overview see Barnett and Lewis [(1994), Chapter 7]. In the case $g \geq 2$, mixture models with outliers have been well known for some time; see again Fraley and Raftery (2002). With the aim of robustifying the *trace* criterion, Cuesta-Albertos, Gordaliza and Matrán (1997) introduced a trimmed version which they called *impartial trimming*: given a trimming level $\alpha \in\, ]0, 1[$, find the subset of the data of size $\lfloor n(1-\alpha) \rfloor$ which is optimal w.r.t. the trace criterion. They also studied its consistency. Later, Garciá-Escudero and Gordaliza (1999) showed robustness properties of the algorithm and, recently, Garciá-Escudero, Gordaliza and Matrán (2003) presented a trimmed $k$-means algorithm for approximating the minimum of the criterion.

We propose first a statistical clustering model with outliers which we call the *spurious-outliers model*. The idea behind it is general enough to allow the derivation of robust clustering criteria with trimming under all kinds of distributional assumptions and cross-cluster constraints. In fact, in the case of normal distributions with equal and spherical covariance matrices, one recovers impartial trimming. Applying the method to equal but general covariance matrices with rejection of $n-r$ elements, the ML-estimator leads us in Section 2 to a robust version of the *pooled determinant* criterion, the *trimmed determinant criterion* (TDC): choose a subset of size $r$ from the $n$ observations and partition it into $g$ clusters so that the pooled SSP matrix has minimum determinant. Not surprisingly, the maximum likelihood estimate of the mean vectors of the different underlying normal distributions are the sample mean vectors of the various clusters, whereas that of the common covariance matrix is the pooled SSP matrix divided by $r$. In the case $r = n$, the TDC simplifies to the classical determinant criterion of cluster analysis.

The number $r$ of regular objects of the model becomes a parameter of the proposed algorithm, the number of retained elements. It turns out that the estimated means and pooled covariance matrix are fairly insensitive to the choice of this number as long as it is not chosen too large. Moreover, we propose in Section 5 a way of estimating $r$ by a method akin to a $\chi^2$ goodness-of-fit test: run the algorithm several times with various values for $r$ and choose the one for which the output best fits the theoretical tail probabilities. This rule may be satisfactory, so much the more as there is no rigorous and unified concept of "outlier," let alone a formal definition; see, for example, Barnett and Lewis (1994) or the introductory discussion in Ritter and Gallegos (1997).

Apart from this parameter, we do not address the question of model finding. Normality of the population distributions, the commonness of the covariance matrix and the number of clusters are assumed as a priori given. This may yield criticism. However, it is straightforward to carry over the method to other *location and scale models*, for example, elliptical distributions with a given radial behavior. But the efficiency of the algorithm depends on that of the ML-estimator of the population parameters and one reason for the popularity of the normal model is the fact that ML-estimation of its parameters essentially reduces to summation. We could also dispense with *commonness* of the covariance matrix. However, the present model should be preferred in situations where each class arises from noisy versions of $g$ prototypes and the noise affects each prototype in the same way. Examples are the classification of phonemes in speech recognition and the chromosome classification problem. In the former case, the prototypes are the phonemes pronounced by a pure speaker and in the latter, they are clean images of the chromosomes of the different biological classes of an organism. In both cases, outliers play an important role; see, for example, Ritter and Gallegos (1997). Moreover, different covariance matrices would require estimation of more parameters and might need more observations than are actually available. Estimation of the *number of clusters* is an important issue that would go beyond the scope of this paper. Just as there is no clear definition of outlier, there is none of "cluster" either. Nevertheless, both are useful concepts. More general distributions, cross-cluster constraints and estimation of the number of clusters in a Bayesian framework will be the subject matter of a forthcoming communication.

Minimizing the TDC requires computing a subset of size $r$ of the $n$ observations and its subsequent partitioning into $g$ clusters (one or several clusters produced by our algorithm may be empty); we call such a partition together with the subset a *configuration*. As in the case of the classical determinant criterion, its computation is infeasible except for small data sets and approximation algorithms which are desirable. In Section 3 we formulate a reduction step that, starting from an arbitrary configuration, yields another configuration with lower or equal determinant of the corresponding pooled SSP matrix; it is based on the Mahalanobis distance. Iterative application of reduction steps until convergence and multistart optimization yield an efficient approximation to the required minimum.

A measure of robustness of an estimator is its *breakdown value* or *breakdown point*: the minimum fraction of bad outliers needed to make it succumb. The *asymptotic breakdown value* is its limit as the number of observations increases to infinity. Estimators with zero asymptotic breakdown value lack robustness. This paper would be incomplete if it did not contain a word about this topic and, in Section 4 we compute the breakdown values of the TDC for the mean vectors and for the pooled SSP matrix. It turns out that the asymptotic breakdown value of the SSP matrix is positive. Mean values, too, are robust w.r.t. data sets that meet a certain condition of cluster separation to be specified in Section 4.2. Both facts plead for robustness of the TDC.

In Section 5 we offer a few simulation studies in order to assess the performance of the proposed algorithm. The error rates obtained compare favorably even to recent studies without outliers; see Coleman and Woodruff (2000).

1.1. *General notation and preliminaries.* Given $g \geq 1$ elements $z_j$, $1 \leq j \leq g$, of a set $F$, $z_1^g$ stands for the $g$-tuple $(z_1, \ldots, z_g) \in F^g$. We write $\mathbf{A} > 0$ ($\mathbf{A} \geq 0$) in order to indicate that a symmetric matrix $\mathbf{A} \in \mathbb{R}^{d \times d}$, $d \geq 1$, is positive (semi-)definite and we denote trace and determinant of $\mathbf{A}$ by $\operatorname{tr} \mathbf{A}$ and $\det \mathbf{A}$, respectively. The $d$-dimensional identity matrix is denoted $\mathbf{I}_d$. The symbol $N_d(\boldsymbol{\mu}, \mathbf{V})$ denotes both the $d$-variate normal distribution with mean vector $\boldsymbol{\mu} \in \mathbb{R}^d$ and covariance matrix $\mathbf{V} \in \mathbb{R}^{d \times d}$ and its Lebesgue density function. $N_d(\boldsymbol{\mu}, \mathbf{V})(\mathbf{x})$ denotes the value of this density at $\mathbf{x} \in \mathbb{R}^d$. The *sum of squares and products matrix* (SSP matrix) $W_E$ of a finite, nonempty subset $E \subseteq \mathbb{R}^d$ with mean $\mathbf{m}_E$ is the matrix

$$
W_E = \sum_{\mathbf{x} \in E} (\mathbf{x} - \mathbf{m}_E)(\mathbf{x} - \mathbf{m}_E)^T. \tag{1}
$$

We next recall some definitions and notation in the theory of cluster analysis. Let $D = \mathbf{x}_1^n = (\mathbf{x}_1, \mathbf{x}_2, \ldots, \mathbf{x}_n)$ be a list of $n$ observations $\in \mathbb{R}^d$. We will often identify the observation $\mathbf{x}_i$ with the corresponding index $i \in 1..n$ and a subset of the list with the corresponding subset of $1..n$. Given a finite set $E$, the notation $\binom{E}{r}$ stands for the system of all subsets of $E$ of size $r$. Let $g \geq 1$ be the number of clusters. If $R$ is a nonempty, finite subset of $1..n$, $\mathcal{C}_g(R)$ denotes the set of all partitions $\mathcal{R} = \{R_1, R_2, \ldots, R_g\}$ of $R$ into $g$ clusters (or subsets), that is, the set of all configurations over the set $R$. Let $\mathcal{R} = \{R_1, R_2, \ldots, R_g\} \in \mathcal{C}_g(R)$ be such a configuration. We will often identify $R_j$ with its index $j$ and $\mathcal{R}$ with the integral interval $1..g$. For a nonempty cluster $R_j$, let $\mathbf{m}_{R_j}$ be its sample mean vector while, for all empty clusters $R_j$, $\mathbf{m}_{R_j}$ is some vector in $\mathbb{R}^d$ given a priori. We write $\mathbf{m}_\mathcal{R} := (\mathbf{m}_{R_1}, \ldots, \mathbf{m}_{R_g})$ and $\mathbf{m}_\mathcal{R}(j) = \mathbf{m}_{R_j}$. The *pooled* or *within-groups sum of squares and products matrix* (SSP matrix) $\mathbf{W}_\mathcal{R}$ of a configuration $\mathcal{R}$ is the sum of the SSP matrices (1) of all clusters,

$$(2) \qquad \mathbf{W}_\mathcal{R} = \sum_{j=1}^g \sum_{\mathbf{x} \in R_j} (\mathbf{x} - \mathbf{m}_{R_j})(\mathbf{x} - \mathbf{m}_{R_j})^T.$$

**2. The spurious-outliers model and its ML-estimator.** Let $r \leq n$ be the assumed number of regular observations. Both the number $g$ of clusters and $r$ are input parameters of the present clustering problem (concerning the choice of $r$, see Section 5, in particular, Table 1).

2.1. *The spurious-outliers model.* This section extends the usual statistical clustering setup; see Mardia, Kent and Bibby [(1979), Section 13.2], combining it with Mathar's [(1981), Section 5.2] outlier model.

Let $(g_\psi)_{\psi \in \Psi}$ be some family of p.d.f.'s on $\mathbb{R}^d$. The parameter set of our statistical model is

$$(3) \qquad \Theta := \left[\bigcup_{R \in \binom{1..n}{r}} \mathcal{C}_g(R)\right] \times \underbrace{(\mathbb{R}^d)^g \times \{\mathbf{V} \in \mathbb{R}^{d \times d} | \mathbf{V} > 0\} \times \Psi^n}_{\Theta_1}.$$

The first factor of $\Theta$ stands for the unknown configuration, the next two for the unknown parameters of the $g$ underlying normally distributed statistical populations which generate the regular observations. Finally, the last factor of $\Theta$ represents the unknown statistical laws that generate the outliers. Let $X_i$, $i \in 1..n$, be $n$ independent, $\mathbb{R}^d$-valued random variables and let their p.d.f.'s conditional on the parameter $\theta = (\mathcal{R}, \boldsymbol{\mu}_1^g, \mathbf{V}, \psi_1^n)$, $\mathcal{R} = \{R_1, \ldots, R_g\} \in \mathcal{C}_g(R)$ be given by

$$f_{X_i}^\theta = \begin{cases} N_d(\boldsymbol{\mu}_j, \mathbf{V}), & i \in R_j, \\ g_{\psi_i}, & i \notin R. \end{cases}$$

The observations $\mathbf{x}_i$ are realizations of these random variables. Since both regular observations and outliers come from Lebesgue-continuous populations, it is natural to assume that the realizations are in general position (any $d+1$ elements are affinely independent). If we additionally require $r > gd$, then the pigeon hole principle ensures that at least one cluster contains $d+1$ or more elements, which implies $\mathbf{W}_{\mathcal{R}} > 0$ for all configurations $\mathcal{R}$.

At the expense of the stronger condition $2r - n > gd$ instead of $r > gd$, the condition "data set in general position" could be relaxed to "$r$ elements of the data set in general position." This, too, would guarantee nonsingularity of all pooled SSP matrices $\mathbf{W}_{\mathcal{R}}$. This modification would allow the $n - r$ outliers to possess any pattern, for example, a regular one. It would, however, exclude examples with many outliers from the beginning, such as in Fraley and Raftery [(2002), Figure 7]. We, therefore, stick to the former conditions. The user may want to screen the data set for affine dependencies in a preprocessing step, at least if dimension is not high. (Otherwise, one may run the algorithm removing affine dependencies as soon as a singular SSP matrix is detected.) Since the regular populations are assumed to be Lebesgue continuous (even normal), such dependencies are the clearest indications of outliers.

The likelihood function of the spurious-outliers model for the data $\mathbf{x}_1^n$ is

$$(4) \qquad L_{\mathbf{x}_1^n}(\mathcal{R}, \boldsymbol{\mu}_1^g, \mathbf{V}, \psi_1^n) = \left[ \prod_{j=1}^g \prod_{i \in R_j} N_d(\boldsymbol{\mu}_j, \mathbf{V})(\mathbf{x}_i) \right] \left[ \prod_{i \notin R} g_{\psi_i}(\mathbf{x}_i) \right].$$

The maximum likelihood estimate of $\mathcal{R}$, $\boldsymbol{\mu}_1^g$, $\mathbf{V}$, $\psi_1^n$ is any element of the parameter set $\Theta$ which maximizes (4). The following proposition shows that the ML-estimator exists and has a simple representation if the outlier model and the data meet a certain condition (5). We were led to this condition and the method of proof of the following proposition by a similar condition appearing in Ritter and Gallegos (2002) in a different statistical context. The corollary following the proposition exposes an outlier model for which this condition is independent of data.

PROPOSITION 2.1 (ML-estimator).

(a) *If, for each subset $T \subseteq 1..n$ of size $n - r$, the function $\prod_{i \in T} g_{\psi_i}(\mathbf{x}_i)$ possesses a maximum w.r.t. $(\psi_i) \in \Psi^T$, then the maximum likelihood estimator of the spurious-outliers model exists.*

*Assume, in addition,*

$$(5) \qquad \underset{\mathcal{R} \in \bigcup_{R \in \binom{1..n}{r}} \mathcal{C}_g(R)}{\arg\min} \det \mathbf{W}_{\mathcal{R}} \subseteq \underset{\mathcal{R} \in \bigcup_{R \in \binom{1..n}{r}} \mathcal{C}_g(R)}{\arg\max} \max_{(\psi)_1^n} \prod_{i \notin R} g_{\psi_i}(\mathbf{x}_i).$$

*(Note that the "arg max" on the right-hand side of the inclusion exists and that the product depends only on the choice of the outliers.) Then:*

(b) *The MLE of the configuration $\mathcal{R}$ (for the given statistical model) is determined by the TDC*

(TDC) $\quad \min\limits_{\mathcal{R} \in \bigcup_{R \in \binom{1..n}{r}} \mathcal{C}_g(R)} \det \mathbf{W}_\mathcal{R} \left( = \min\limits_{R \in \binom{1..n}{r}} \min\limits_{\mathcal{R} \in \mathcal{C}_g(R)} \det \mathbf{W}_\mathcal{R} \right).$

(c) *If we denote it by $\mathcal{R}^\star$ (if it is not unique, choose one), then the MLE of $(\boldsymbol{\mu}_1, \ldots, \boldsymbol{\mu}_g)$ is $(\mathbf{m}_{\mathcal{R}^\star}(1), \ldots, \mathbf{m}_{\mathcal{R}^\star}(g))$ and that of the common covariance matrix $\mathbf{V}$ is $\frac{1}{r}\mathbf{W}_{\mathcal{R}^\star}$.*

PROOF. By the assumptions and the discussion at the beginning of this section, we have $\mathbf{W}_\mathcal{R} > 0$ for any $\mathcal{R} \in \bigcup_{R \in \binom{1..n}{r}} \mathcal{C}_g(R)$. Therefore, classical normal distribution theory [see Mardia, Kent and Bibby (1979), pages 103–105] shows

$$\max_{\mathbf{V}} \max_{\boldsymbol{\mu}_1^g} \prod_{j=1}^{g} \prod_{i \in R_j} N_d(\boldsymbol{\mu}_j, \mathbf{V})(\mathbf{x}_i)$$

$$= \max_{\mathbf{V}} (\det 2\pi \mathbf{V})^{-r/2} \prod_{j=1}^{g} \exp -\tfrac{1}{2} \operatorname{tr}\left( \mathbf{V}^{-1} \sum_{i \in R_j} (\mathbf{x}_i - \mathbf{m}_\mathcal{R}(j))(\mathbf{x}_i - \mathbf{m}_\mathcal{R}(j))^T \right)$$

(6)
$$= \max_{\mathbf{V}} (\det 2\pi \mathbf{V})^{-r/2} \exp -\tfrac{1}{2} \operatorname{tr}(\mathbf{V}^{-1} \mathbf{W}_\mathcal{R})$$

$$= K (\det 2\pi \mathbf{W}_\mathcal{R})^{-r/2},$$

where $K$ is a constant independent of $\mathcal{R}$. The last equality is a direct application of Mardia, Kent and Bibby [(1979), Theorem 4.2.1, page 104]. Finiteness of $\bigcup_{R \in \binom{1..n}{r}} \mathcal{C}_g(R)$, (6), (4) and the hypothesis of (a) together imply

(7)
$$\max_{\mathcal{R}} \left[ \max_{\mathbf{V}} \max_{\boldsymbol{\mu}_1^g} \prod_{j=1}^{g} \prod_{i \in R_j} N_d(\boldsymbol{\mu}_j, \mathbf{V})(\mathbf{x}_i) \max_{\psi_1^n} \prod_{i \notin R} g_{\psi_i}(\mathbf{x}_i) \right]$$
$$= \max_{(\mathcal{R}, \boldsymbol{\mu}_1^g, \mathbf{V}, \psi_1^n)} L_{\mathbf{x}_1^n}(\mathcal{R}, \boldsymbol{\mu}_1^g, \mathbf{V}, \psi_1^n),$$

which proves part (a).

Now, by hypothesis (5), any configuration $\mathcal{R}$ that maximizes the first factor in brackets in (7) also maximizes the second one. The maximum is, therefore, attained at the minimum of $\det \mathbf{W}_\mathcal{R}$ over all configurations $\mathcal{R}$. □

In order to minimize the criterion (TDC), one has to determine a *subset* $R \subseteq 1..n$ of $r$ elements, the set of estimated *regular* observations and a *partition* $\mathcal{R}$ of it in $g$ clusters such that $\mathbf{W}_\mathcal{R}$ has minimum determinant. Since

$$(6) = \prod_{j=1}^{g} \prod_{i \in R_j} N_d\left(\mathbf{m}_\mathcal{R}(j), \frac{1}{r}\mathbf{W}_\mathcal{R}\right)(\mathbf{x}_i),$$

this may be restated as follows: in order to find $\mathcal{R}$, compute the subset of $r$ observations which is optimally explained by $g$ normal populations with a common sample covariance.

Condition (5) is not very restrictive. It is satisfied if $\max_\psi g_\psi(\mathbf{x}_i)$ exists for all $i$ and does not depend on $i$. Two special, opposite cases are the following.

COROLLARY 2.1. *The conclusions of Proposition* 2.1 *hold under each of the following two conditions:*

(a) $\Psi = \mathbb{R}^d$ *and* $g_\psi(\mathbf{x}) = g(\mathbf{x} - \psi)$ *is a location model with a density function* $g$ *having a maximum* [*e.g.*, $g = N_d(0, \mathbf{I}_d)$ *or* $g = \frac{1}{\lambda^d(B)} \cdot \mathbf{1}_B$, *where* $B$ *is some region about the origin*].
(b) $\Psi$ *is singleton and* $g_\psi$ *is constant on a region that contains all data.*

PROOF. In case (a), let $M = \max g$. Clearly, from the assumption,

$$\max_{\psi_1^n} \prod_{i \notin R} g_{\psi_i}(\mathbf{x}_i) = \prod_{i \notin R} \max_{\psi_i} g(\mathbf{x}_i - \psi_i) = M^{n-r}$$

does not depend on $\mathcal{R}$ and Condition (5) of Proposition 2.1 is satisfied. The same is obviously true in case (b). □

The corollary shows that the parameter set $\Psi$ may be of any size. If it contains just one element, then all outliers belong to *one* population. The model also allows that each outlier comes from its *own* population. This particularity justifies the adjective "spurious."

Criterion (TDC) is optimal (in the maximum likelihood sense) for the spurious-outliers model 2.1 if condition (5) is satisfied, in particular, if the outliers are generated from Corollary 2.1(a) or (b). Therefore, it performs well whenever the data set $\mathbf{x}_1^n$ is a realization of random variables meeting this condition. However, it is also a plausible descriptive tool per se. We formulate the case of one cluster (robust parameter estimation) as a second corollary of Proposition 2.1. A similar statement for normally distributed outliers appears in Pesch (2000) where Mathar's outlier model is already used.

COROLLARY 2.2. *Assume that the data* $\mathbf{x}_1^n$ *satisfies* (5) [*e.g., that the family* $(g_\psi)_{\psi \in \mathbb{R}^d}$ *is a location model*]. *Then Rousseeuw's MCD is the maximum likelihood estimator for the spurious-outliers model* 2.1 *with* $g = 1$.

REMARKS. (a) *Classical determinant criterion*. For $r = n$ (the pure clustering situation), criterion (TDC) reduces to the classical determinant criterion of cluster analysis [Friedman and Rubin (1967) and Scott and Symons (1971)].

(b) *Mixture model vs. classification model.* There exist two approaches to classical (nonrobust) model-based clustering: mixture modeling and the classification approach. It is well known that both are related; see Mardia, Kent and Bibby [(1979), remark (4), page 365] and Fraley and Raftery (2002). The same can be said about mixture modeling with outliers and the present robust classification model 2.1. Indeed, define the parameter set of the mixture model with outliers as

$$\Theta_M := \left\{ (\pi_1, \ldots, \pi_g) \in [0,1]^g \Big| \sum_{j=1}^{g} \pi_j = \frac{r}{n} \right\} \times \Theta_1,$$

where $(\pi_1, \ldots, \pi_g)$ are the mixing parameters and $\Theta_1$ is defined in (3). Let $Y_1, \ldots, Y_n$ be i.i.d., $0..g$-valued random variables whose common distribution under the parameter $\theta = (\pi_1^g, \boldsymbol{\mu}_1^g, \mathbf{V}, \psi_1^n) \in \Theta_M$ is given by $P^\theta[Y_1 = j] = \pi_j$, $j \in 1..g$, and $P^\theta[Y_1 = 0] = 1 - \frac{r}{n}$. Furthermore, let $X_{i,j}$, $i \in 1..n$, $j \in 0..g$, be $n \cdot (g+1)$ independent, $\mathbb{R}^d$-valued random variables, independent of $Y_1^n$ and with p.d.f. conditional on $\theta$,

$$f_{X_{i,j}}^\theta = \begin{cases} N_d(\boldsymbol{\mu}_j, \mathbf{V}), & j \in 1..g, \\ g_{\psi_i}, & j = 0. \end{cases}$$

The formula of total probability shows that the $\mathbb{R}^d$-valued random variables $X_{1,Y_1}, \ldots, X_{n,Y_n}$ are independent and that the conditional p.d.f. of $X_{i,Y_i}$, $i \in 1..n$, is given by the mixture

$$(8) \qquad f_{X_{i,Y_i}}^\theta(\mathbf{x}) = \sum_{j=1}^{g} \pi_j N_d(\boldsymbol{\mu}_j, \mathbf{V})(\mathbf{x}) + \left(1 - \frac{r}{n}\right) g_{\psi_i}(\mathbf{x})$$

with mixing parameters $\pi_j$. The data $\mathbf{x}_1^n$ is now a realization of the random variables $X_{1,Y_1}, \ldots, X_{n,Y_n}$. Note that $1 - \frac{r}{n}$ is the prior probability of occurrence of an outlier. Hence, in this model the configuration $\mathcal{R}$ is an unobservable random variable. In the special case of Corollary 2.1(b), one finds the mixture model appearing in equation (11) of Fraley and Raftery (2002).

The aim of the statistical clustering model 2.1 and the mixture model is (robust) clustering and estimation of the means of all subpopulations and of their common covariance matrix. Whereas, in doing so, the former estimates (besides these parameters) the optimal configuration $\mathcal{R}^*$, the latter estimates the probabilities of the observations to come from the different clusters. In this sense, both models pursue the same aim. In the clustering model, the prior information of the existence of $n - r$ outliers is expressed by the constraint $\#\complement R = n - r$. The mixture model describes this fact by setting the probability of occurrence of a contamination to $1 - \frac{r}{n}$. Furthermore, once the ML-estimates $(\hat{\pi}_1^g, \hat{\mu}_1^g, \widehat{\mathbf{V}})$ of the mixture (8) are known, we can regard each distribution $N_d(\hat{\mu}_j, \widehat{\mathbf{V}})$ as indicating a separate group, and individuals are then assigned to clusters by Bayes' allocation rule: assign the $i$th observation, $i \in 1..n$, to the class $j$ that maximizes the posterior density $\hat{\pi}_j N_d(\hat{\mu}_j, \widehat{\mathbf{V}})(\mathbf{x}_i)$. The (estimated) set of regular observations $R$ consists of the $r$ observations with the $r$ largest maxima. The corresponding optimal partition is defined by the class assignments of the elements of $R$. Similarly, once the optimal configuration $\mathcal{R}^*$ of the clustering model is known, we can estimate the mixing parameters by the cluster sizes divided by $n$.

(c) *Unequal cluster sizes.* Being an ML-estimator, the pooled determinant criterion can be interpreted as a maximum a posteriori estimator for mixtures with *equal* mixing parameters. Therefore, it favors equal cluster sizes, although it can deal with small deviations from the ideal situation. For *unequal* mixing parameters, an entropy correction has to be added to the criterion; see Symons (1981). The same remark applies to the trimmed version. We will deal with this topic (and the related question of the number of clusters) in another communication.

**3. An efficient approximation algorithm.** Minimizing the trimmed determinant criterion requires the computation of a subset of size $r$ out of the $n$ observations and its subsequent partitioning into $g$ clusters. This task is infeasible, except for small data sets, and an efficient approximation algorithm is desirable. In this section we develop such a procedure. It is iterative and adapts the idea of minimal distance partition, now classical in cluster analysis [see Schroeder (1976) for a general version], to the case with outliers. In the classical case (*without* trimming), one reduces the sum of the squared Mahalanobis distances w.r.t. $\mathbf{W}_\mathcal{R}$ for a fixed "current" configuration $\mathcal{R}$ by reassigning *single observations* to cluster centers with smaller Mahalanobis distances w.r.t. $\mathbf{W}_\mathcal{R}$. Moreover, one shows that this reduction also reduces the determinant of the SSP matrix. We prove below that the same idea can be applied also to the trimmed determinant criterion; this extension is, however, not straightforward. The following theorem gives rise to the basic reduction step of our algorithm.

THEOREM 3.1. *Let $\mathcal{R}$ and $\mathcal{R}_{\text{new}}$ be two configurations over $r$-element subsets $R$, $R_{\text{new}} \subseteq 1..n$, respectively, such that*

$$\sum_{j=1}^{g} \sum_{i \in R_{\text{new},j}} (\mathbf{x}_i - \mathbf{m}_\mathcal{R}(j))^T \mathbf{W}_\mathcal{R}^{-1} (\mathbf{x}_i - \mathbf{m}_\mathcal{R}(j))$$

(9)

$$\leq \sum_{j=1}^{g} \sum_{i \in R_j} (\mathbf{x}_i - \mathbf{m}_{\mathcal{R}}(j))^T \mathbf{W}_{\mathcal{R}}^{-1} (\mathbf{x}_i - \mathbf{m}_{\mathcal{R}}(j)).$$

(a) *We have* $\det \mathbf{W}_{\mathcal{R}_{\text{new}}} \leq \det \mathbf{W}_{\mathcal{R}}$ *with equality if and only if* $\mathbf{W}_{\mathcal{R}_{\text{new}}} = \mathbf{W}_{\mathcal{R}}$.

(b) *Let us put* $\mathbf{m}_{\mathcal{R}_{\text{new}}}(j) := \mathbf{m}_{\mathcal{R}}(j)$ *for all* $j \in 1..g$ *such that* $R_{\text{new},j} = \varnothing$. *If there is equality in* (a), *then we have also* $\mathbf{m}_{\mathcal{R}_{\text{new}}} = \mathbf{m}_{\mathcal{R}}$.

PROOF. At the beginning of the proof of Proposition 2.1, we have already used the fact that the determinant of $\mathbf{W}_{\mathcal{R}}$ is a constant multiple of a negative power of the product

$$\prod_{j=1}^{g} \prod_{i \in R_j} N_d\left(\mathbf{m}_{\mathcal{R}}(j), \frac{1}{r}\mathbf{W}_{\mathcal{R}}\right)(\mathbf{x}_i).$$

Claim (a) will, therefore, follow if we prove

$$\prod_{j=1}^{g} \prod_{i \in R_{\text{new},j}} N_d\left(\mathbf{m}_{\mathcal{R}_{\text{new}}}(j), \frac{1}{r}\mathbf{W}_{\mathcal{R}_{\text{new}}}\right)(\mathbf{x}_i) \geq \prod_{j=1}^{g} N_d\left(\mathbf{m}_{\mathcal{R}}(j), \frac{1}{r}\mathbf{W}_{\mathcal{R}}\right)(\mathbf{x}_i).$$

Now, passing to likelihoods, we have

$$\prod_{j=1}^{g} \prod_{i \in R_{\text{new},j}} L_{\mathbf{x}_i}\left(\mathbf{m}_{\mathcal{R}_{\text{new}}}(j), \frac{1}{r}\mathbf{W}_{\mathcal{R}_{\text{new}}}\right)$$

$$\geq \prod_{j=1}^{g} \prod_{i \in R_{\text{new},j}} L_{\mathbf{x}_i}\left(\mathbf{m}_{\mathcal{R}}(j), \frac{1}{r}\mathbf{W}_{\mathcal{R}}\right)$$

(10)
$$= \exp \sum_{j} \sum_{i \in R_{\text{new},j}} l_{\mathbf{x}_i}\left(\mathbf{m}_{\mathcal{R}}(j), \frac{1}{r}\mathbf{W}_{\mathcal{R}}\right)$$

$$\geq \exp \sum_{j} \sum_{i \in R_j} l_{\mathbf{x}_i}\left(\mathbf{m}_{\mathcal{R}}(j), \frac{1}{r}\mathbf{W}_{\mathcal{R}}\right)$$

$$= \prod_{j=1}^{g} \prod_{i \in R_j} L_{\mathbf{x}_i}\left(\mathbf{m}_{\mathcal{R}}(j), \frac{1}{r}\mathbf{W}_{\mathcal{R}}\right).$$

In this chain the first inequality follows from ML-estimation and the second is just the assumption. This proves the first part of (a).

If the two determinants are equal, so are both ends of the above chain. Equality in (10) follows and uniqueness of the ML-estimator implies $\mathbf{W}_{\mathcal{R}_{\text{new}}} = \mathbf{W}_{\mathcal{R}}$ and $\mathbf{m}_{\mathcal{R}_{\text{new}}}(j) = \mathbf{m}_{\mathcal{R}}(j)$ for all $j \in 1..g$ such that $R_{\text{new},j} \neq \varnothing$. This concludes the proofs of parts (a) and (b). $\square$

Let $\mathcal{R}$ be any configuration. With the squared Mahalanobis distances

(11) $\quad d_{\mathcal{R}}(i,j)^2 := (\mathbf{x}_i - \mathbf{m}_{\mathcal{R}}(j))^T \mathbf{W}_{\mathcal{R}}^{-1}(\mathbf{x}_i - \mathbf{m}_{\mathcal{R}}(j)), \qquad i \in 1..n, j \in 1..g.$

Inequality (9) may be rewritten as

(12) $$\sum_{j=1}^{g} \sum_{i \in R_{\text{new},j}} d_{\mathcal{R}}(i,j)^2 \leq \sum_{j=1}^{g} \sum_{i \in R_j} d_{\mathcal{R}}(i,j)^2.$$

Theorem 3.1 is the basis of the following building block for our algorithm: starting from a configuration $\mathcal{R}$, look for another configuration $\mathcal{R}_{\text{new}}$ such that the corresponding sum of distance squares w.r.t. *the given* $\mathcal{R}$ is smaller than the current one; see (12). The theorem assures that this new configuration is a better approximation to the minimum of the TDC.

Plainly, a configuration $\mathcal{R}_{\text{new}}$ that minimizes the sum

(13) $$\sum_{j=1}^{g} \sum_{i \in P_j} d_{\mathcal{R}}(i,j)^2$$

over all configurations $\{P_1, \ldots, P_g\} \in \mathcal{C}_g(P)$, $P \in \binom{1..n}{r}$, satisfies (12). Fortunately, computing this minimum is very simple: it is sufficient to assign each observation $i$ to a cluster $j \in 1..g$ which minimizes the distance square $d_{\mathcal{R}}(i,j)^2$ with respect to the fixed configuration $\mathcal{R}$. Let us call each cluster $j \in 1..g$ that minimizes $d_{\mathcal{R}}(i,j)^2$ *optimal* for the $i$th observation, $i \in 1..n$, with respect to the given partition $\mathcal{R}$. Since we must restrict our choice to $r$ observations, the optimal ones are those with the $r$ smallest distances to their optimal clusters. These ideas are made precise in the following corollary of Theorem 3.1.

COROLLARY 3.1. *Let $\mathcal{R}$ be a configuration and let $R_{\text{new}}$ be a subset of $1..n$ consisting of $r$ observations with the smallest Mahalanobis distances to their optimal clusters w.r.t. $\mathcal{R}$ (in general, $R_{\text{new}}$ is unique). Let $\mathcal{R}_{\text{new}}$ be the partition of $R_{\text{new}}$ obtained by assigning each $i \in R_{\text{new}}$ to its optimal cluster w.r.t. $\mathcal{R}$. Then:*

(a) $\mathcal{R}_{\text{new}}$ *minimizes the objective function* (13) *over the set of all configurations* $\{P_1, \ldots, P_g\} \in \mathcal{C}_g(P)$, $P \in \binom{1..n}{r}$.

(b) *The conclusions of Theorem 3.1 hold.*

In the case of one class, $g = 1$, Corollary 3.1(b) is the basis of Rousseeuw and Van Driessen's C-step [Rousseeuw and Van Driessen (1999), Theorem 1]. For $n = r$, (the noncontaminated situation) Corollary 3.1(b) reduces to Späth [(1985), Theorem 3.5]. We next formulate the reduction step in algorithmic terms.

3.1. *The reduction step.* *Input*: A configuration $\mathcal{R}$ together with its mean vectors $\mathbf{m}_\mathcal{R}$ and its SSP matrix $\mathbf{W}_\mathcal{R}$;

*Output*: A configuration $\mathcal{R}_{\text{new}}$ such that $\det \mathbf{W}_{\mathcal{R}_{\text{new}}} \leq \det \mathbf{W}_\mathcal{R}$.

(i) Compute the distance squares $d_\mathcal{R}(i,j)^2$, $i \in 1..n$, $j \in 1..g$, defined in (11).

(ii) For each $i \in 1..n$, determine an optimal cluster $j_i \in \mathcal{R}$, that is, $j_i \in \arg\min_{j \in 1..g} d_\mathcal{R}(i,j)^2$.

(iii) Determine a permutation $\kappa : 1..n \to 1..n$ that satisfies

$$(14) \quad d_\mathcal{R}(\kappa(1), j_{\kappa(1)})^2 \leq d_\mathcal{R}(\kappa(2), j_{\kappa(2)})^2 \leq \cdots \leq d_\mathcal{R}(\kappa(n), j_{\kappa(n)})^2.$$

(iv) Put $R_{\text{new}} = \{\kappa(1), \ldots, \kappa(r)\}$ and, for each $j \in 1..g$, put $R_{\text{new},j} = \{i \in 1..r \mid j_{\kappa(i)} = j\}$. Finally, let $\mathcal{R}_{\text{new}} := \{R_{\text{new},1}, \ldots, R_{\text{new},g}\}$.

3.2. *Iteration and discussion.* Now, starting from an initial configuration $\mathcal{R}_0$ and iterating reduction steps, we obtain a sequence of configurations $(\mathcal{R}_k)_{k \geq 0}$ such that $\det \mathbf{W}_{\mathcal{R}_{k+1}} \leq \det \mathbf{W}_{\mathcal{R}_k}$ for all $k$. Since there are only a finite number of configurations, this iterative process must become stationary after a finite number of steps, say $L$, with $\det \mathbf{W}_{\mathcal{R}_{L+1}} = \det \mathbf{W}_{\mathcal{R}_L}$ ($> 0$). By Corollary 3.1, we have $\mathbf{W}_{\mathcal{R}_{L+1}} = \mathbf{W}_{\mathcal{R}_L}$ and $\mathbf{m}_{\mathcal{R}_{L+1}} = \mathbf{m}_{\mathcal{R}_L}$. Therefore, $d_{\mathcal{R}_L}(i,j) = d_{\mathcal{R}_{L+1}}(i,j)$, $i \in 1..n$, $j \in 1..g$, and, if $\mathcal{R}_L$ is unique, a further reduction step yields a configuration with sum (13) and the TDC unchanged. If $\mathcal{R}_L$ is not unique, then a further step may improve the TDC, but not (13). (An example of nonuniqueness in the complex plane is $\mathbf{x}_k = e^{i\pi k/4}$, $k \in 1..8$, $r = 4$, $g = 1$ and $\mathcal{R}_0 = \{\mathbf{x}_1, \mathbf{x}_3, \mathbf{x}_5, \mathbf{x}_7\}$.)

The configuration $\mathcal{R}_L$ is *one* approximation to the minimum. Now, multistart optimization is applied to the foregoing iterative process; the limit configuration with the least value of the determinant of the corresponding SSP matrix is the final approximation to the minimum.

If a configuration $\mathcal{R}$ is a *global minimum* of the TDC, then a reduction step with input $\mathbf{m}_\mathcal{R}$ and $\mathbf{W}_\mathcal{R}$ yields an equivalent configuration.

3.3. *The initial configuration.* We indicate two methods for generating random initial configurations. Both are natural extensions of the ones proposed by Rousseeuw and Van Driessen [(1999), Section 4.1]:

(a) Draw a random configuration consisting of nonempty clusters.

(b) Choose at random a subset of $1..n$ with at least $gd + 1$ elements. Construct a random partition $\mathcal{D}$ of the subset in $g$ clusters and compute its mean vectors $\mathbf{m}_\mathcal{D}$ and its SSP matrix $\mathbf{W}_\mathcal{D}$. Iterate a reduction step to obtain an initial configuration $\mathcal{R}_0$.

We conclude this section with a result concerning some geometrical properties of any limit configuration of our algorithm. This result extends Corollary 1 of Rousseeuw and Van Driessen (1999) to robust clustering and the well-known geometric separation property of discriminant analysis, see Mardia, Kent and Bibby [(1979), Theorem 11.2.1], to clustering in the presence of outliers.

COROLLARY 3.2. *Let $\mathcal{R} = \{R_1, \ldots, R_g\}$ be a limit configuration of the reduction step iteration, for example, the optimal configuration.*

(a) *Each nonempty cluster $R_j$, $j \in 1..g$, is separated from the estimated set $\complement \bigcup_{j=1}^{g} R_j$ of outliers by an ellipsoid.*

(b) *Two different, nonempty clusters $R_j$ and $R_l$ are separated by the hyperplane $h_{jl}(\mathbf{x}_i) = 0$, where $h_{jl} : \mathbb{R}^d \to \mathbb{R}$ is the linear form*

$$h_{jl}(\mathbf{y}) := 2[\mathbf{y} - \tfrac{1}{2}(\mathbf{m}_{\mathcal{R}}(j) + \mathbf{m}_{\mathcal{R}}(l))]^T \mathbf{W}_{\mathcal{R}}^{-1}(\mathbf{m}_{\mathcal{R}}(j) - \mathbf{m}_{\mathcal{R}}(l)), \qquad \mathbf{y} \in \mathbb{R}^d.$$

*The observations $i$ in cluster $R_j$ are those satisfying $h_{jl}(\mathbf{x}_i) \geq 0$.*

PROOF. The application of a reduction step to the configuration $\mathcal{R}$ yields $\mathcal{R}$ itself as possible output. Thus, the set of regular observations $R = \bigcup_{j=1}^{g} R_j$ may be written as $\{\kappa(1), \kappa(2), \ldots, \kappa(r)\}$, where $\kappa : 1..n \to 1..n$ is a permutation satisfying (14), whereas the set of outliers is given by $\{\kappa(r+1), \kappa(r+2), \ldots, \kappa(n)\}$. In order to prove part (a), let $j \in 1..g$ such that $R_j \neq \varnothing$. All observations $i \in R_j$ satisfy $d_{\mathcal{R}}(i,j)^2 \leq \max_{1 \leq m \leq r} d_{\mathcal{R}}(\kappa(m), j)^2 =: K_j$, whereas all $i \notin R$ satisfy $d_{\mathcal{R}}(i, j_i)^2 \geq K_j$ [even $d_{\mathcal{R}}(i, j_i)^2 \geq d_{\mathcal{R}}(\kappa(r), j_{\kappa(r)})^2$]. The ellipsoid

$$E_j = \{\mathbf{x} \in \mathbb{R}^d | (\mathbf{x} - \mathbf{m}_{\mathcal{R}}(j))^T \mathbf{W}_{\mathcal{R}}^{-1}(\mathbf{x} - \mathbf{m}_{\mathcal{R}}(j)) \leq K_j\}$$

contains $R_j$ and $\complement R$ is contained in the closure of $\complement E_j$.

For two observations $i_1$ and $i_2$ in the $j$th and in the $l$th cluster, respectively, we have $d_{\mathcal{R}}(i_1, j)^2 \leq d_{\mathcal{R}}(i_1, l)^2$ and $d_{\mathcal{R}}(i_2, l)^2 \leq d_{\mathcal{R}}(i_2, j)^2$. Part (b) now follows from standard matrix operations. □

## 4. The breakdown values.

4.1. *Preliminaries.* Besides the asymptotic breakdown value of an estimator Hampel (1968, 1971), there exists also a finite-sample version, Hodges (1967) and Donoho and Huber (1983). Loosely speaking, the latter measures the minimum fraction of bad outliers that can *completely* spoil the estimate. More precisely, let $\Theta$ be a locally compact parameter space, for example, the intersection of an open and a closed subset of some Euclidean space and consider an estimator $\delta : \mathcal{A} \to \Theta$. Here, $\mathcal{A} \subseteq \mathbb{R}^{n \cdot d}$ is the system of all data sets *admissible* for $\delta$ ("in general position" in our case). We say that $(\mathbf{x}'_1, \ldots, \mathbf{x}'_n) \in \mathcal{A}$ is an $m$-modification, $m \leq n$, of a data set $(\mathbf{x}_1, \ldots, \mathbf{x}_n) \in \mathcal{A}$ if it arises from $(\mathbf{x}_1, \ldots, \mathbf{x}_n)$ by replacing $m$ observations $\mathbf{x}_i$ with arbitrary elements $\mathbf{x}'_i \in \mathbb{R}^d$ in an admissible way. An estimator $\delta : \mathbb{R}^{n \cdot d} \to \Theta$ *breaks* *down* under $m$ replacements with a data set $(\mathbf{x}_1,\ldots,\mathbf{x}_n) \in \mathcal{A}$ if the set of estimates

$$\{\delta(\mathbf{x}'_1,\ldots,\mathbf{x}'_n) | (\mathbf{x}'_1,\ldots,\mathbf{x}'_n) \text{ is an } m\text{-modification of } (\mathbf{x}_1,\ldots,\mathbf{x}_n)\} \subseteq \Theta$$

is not relatively compact in $\Theta$. The *individual breakdown point* for $\mathbf{x}_1^n$ is defined as

$$\beta(\delta, \mathbf{x}_1^n) = \min_{1 \leq m \leq n}\left\{\frac{m}{n}\Big| \{\delta(M)|M\} \text{ is not relatively compact in } \Theta\right\};$$

here $M$ runs over all $m$-modifications of $\mathbf{x}_1^n$. It is the minimum fraction of replacements in $\mathbf{x}_1^n$ that may cause $\delta$ to break down.

Depending on a specific data set, this is not a useful notion per se. Therefore, Donoho and Huber define a value that we call the *universal breakdown value* $\beta(\delta)$ of $\delta$: it is the minimum relative amount of replacements that causes $\delta$ to break down with *some* data set $\mathbf{x}_1^n \in \mathcal{A}$:

(15) $$\beta(\delta) = \min_{\mathbf{x}_1^n \in \mathcal{A}} \beta(\delta, \mathbf{x}_1^n).$$

According to this definition, the estimator breaks down at the first integer $m$ for which *there exists* some $\mathbf{x}_1^n$ such that the estimate becomes arbitrarily bad for some suitable modification $M$.

The universal breakdown value is the minimal individual one; it depends on the estimator and its parameters alone, not on data. It is pessimistic in considering the *worst case* and modifications of this notion are conceivable. One may argue that the existence of a *single* data set $\mathbf{x}_1^n$ and possibly very *special* bad modifications $M$ may not suffice to indicate lack of robustness of an estimator. A more relaxed definition would require the criterion to break down for *sufficiently many* data sets $\mathbf{x}_1^n$. We will introduce such a modification in Section 4.2. Another less stringent definition would require all components of the estimate to break down [such as all means in Definition 4.1(a)].

The present task is, among other things, estimating mean vectors and an SSP matrix by means of the TDC. In the first case, $\Theta = \mathbb{R}^{g \cdot d}$ and in the second, $\Theta$ is the set of all positive-definite $d$ by $d$ matrices, an open subset of $\binom{d+1}{2}$-dimensional Euclidean space.

Our definitions and analyses need the following facts. If $\mathbf{A} \leq \mathbf{B}$, then $\operatorname{tr} \mathbf{A} \leq \operatorname{tr} \mathbf{B}$ by linearity and $\det \mathbf{A} \leq \det \mathbf{B}$ (see Lemma A.2). Let $\lambda_{\min}(\mathbf{A})$ ($\lambda_{\max}(\mathbf{A})$) be the least (largest) eigenvalue of a matrix $\mathbf{A} \geq 0$. Then

$$\lambda_{\min}(\mathbf{A}) = \min_{\|x\|=1} x^T \mathbf{A} x \leq \min_{\|x\|=1} x^T \mathbf{B} x = \lambda_{\min}(\mathbf{B})$$

and, similarly,

$$\lambda_{\max}(\mathbf{A}) \leq \lambda_{\max}(\mathbf{B}).$$

Let $E \subseteq F$ be nonempty, finite subsets of $\mathbb{R}^d$ and let $\mathbf{m}_E$ and $\mathbf{m}_F$ be their means. Then their SSP matrices $W_E$ and $W_F$ [in the simple sense (1)] satisfy $W_E = W_E(\mathbf{m}_E) \leq W_E(\mathbf{m}_F) \leq W_F(\mathbf{m}_F) = W_F$. Hence, also $\operatorname{tr} W_E \leq \operatorname{tr} W_F$, $\det W_E \leq \det W_F$, $\lambda_{\min}(W_E) \leq \lambda_{\min}(W_F)$ and $\lambda_{\max}(W_E) \leq \lambda_{\max}(W_F)$.

In the present situation, corruption of the estimates is reflected by an arbitrarily large value of at least one sample mean or by an arbitrarily large or small eigenvalue of the pooled SSP matrix of the optimal configuration; see Donoho and Huber (1983). Transferring definition (15) to the present situation, we obtain the following.

DEFINITION 4.1.　Let $n, r, g, d$ be such that $n \geq r \geq gd + 1$ as before. Given a data set $M \subseteq \mathbb{R}^{n \cdot d}$ of observations in general position, let $\mathcal{M}^\star$ denote its optimal configuration w.r.t. the TDC.

(a) The *universal breakdown value* of the TDC for the *mean vectors* is

$$\beta_{\text{mean}}(n, r, g) = \min_{\mathbf{x}_1^n} \min_{1 \leq m \leq n} \left\{ \frac{m}{n} \Big| \sup_M \max_{j \in 1..g} \|\mathbf{m}_{\mathcal{M}^\star}(j)\| = \infty \right\};$$

here $\mathbf{x}_1^n$ runs over all data sets in general position and $M$ over all $m$-modifications of $\mathbf{x}_1^n$ in general position.

(b) The *universal breakdown value* of the TDC for the SSP matrix is

$$\beta_{\text{SSP}}(n, r, g) = \min_{\mathbf{x}_1^n} \min_{1 \leq m \leq n} \left\{ \frac{m}{n} \Big| \sup_M \max\left(\lambda_{\max}(\mathbf{W}_{\mathcal{M}^\star}), \frac{1}{\lambda_{\min}(\mathbf{W}_{\mathcal{M}^\star})}\right) = \infty \right\},$$

where $\mathbf{x}_1^n$ and $M$ are as stated in (a).

We first show that, in general, the *universal* breakdown value of the TDC w.r.t. the mean vectors is low. We need two lemmas; the first—a geometrical interpretation of the SSP matrix—is of general interest and the second is combinatorial and of a technical nature. The parallelepiped spanned by $k+1$ points $\mathbf{y}_0, \ldots, \mathbf{y}_k \in \mathbb{R}^d$, $k \leq d$, is the subset

$$P(\mathbf{y}_0, \ldots, \mathbf{y}_k) = \left\{ \mathbf{y}_0 + \sum_{i=1}^k \lambda_i(\mathbf{y}_i - \mathbf{y}_0) \Big| 0 \leq \lambda_i \leq 1 \right\} \subseteq \mathbb{R}^d.$$

Its $d$-dimensional volume is independent of the order of the points $\mathbf{y}_i$.

LEMMA 4.1.　*Let $E = \{\mathbf{x}_0, \ldots, \mathbf{x}_d\} \subseteq \mathbb{R}^d$. We have the equality $\det W_E = \frac{1}{d+1} \operatorname{volume}^2 P(E)$.*

PROOF.　Let $\mathbf{m} = \frac{1}{d+1} \sum_{i=0}^d \mathbf{x}_i$. Putting

$$A = \begin{pmatrix} 1 & \cdots & 1 \\ \mathbf{x}_0 - \mathbf{m} & \cdots & \mathbf{x}_d - \mathbf{m} \end{pmatrix} \in \mathbb{R}^{(d+1) \times (d+1)},$$

we obtain the claim from

$$\text{volume}^2 P(\mathbf{x}_0, \ldots, \mathbf{x}_d)$$
$$= \det^2(\mathbf{x}_1 - \mathbf{x}_0, \ldots, \mathbf{x}_d - \mathbf{x}_0)$$
$$= \det^2 A = \det A \det A^T = \det \begin{pmatrix} 1 & \cdots & 1 \\ \mathbf{x}_0 - \mathbf{m} & \cdots & \mathbf{x}_d - \mathbf{m} \end{pmatrix} \begin{pmatrix} 1 & (\mathbf{x}_0 - \mathbf{m})^T \\ \vdots & \vdots \\ 1 & (\mathbf{x}_d - \mathbf{m})^T \end{pmatrix}$$
$$= \det \begin{pmatrix} d+1 & 0 \\ 0 & W_E \end{pmatrix} = (d+1) \det W_E. \qquad \square$$

LEMMA 4.2. *Let $g \geq 2$, $p \geq 2$, $q \geq g-2$ and $r = p+g$ be natural numbers and let*

$$F = \{x_1, \ldots, x_p\} \cup \{y_1, y_2\} \cup \{z_1, \ldots, z_q\}$$

*with pairwise disjoint elements $x_i$, $y_k$ and $z_l$. Any partition $\mathcal{R}$ of a subset of $F$ of size $r$ in $g$ classes is either of the form*

$$\mathcal{R}^\star = \{\{x_1, \ldots, x_p\}, \{y_1, y_2\}, g-2 \text{ one-point classes from the } z_l\text{'s}\}$$

*or possesses a class $C$ which contains some pair $\{x_i, y_k\}$ or some pair $\{z_l, u\}$, $u \neq z_l$.*

PROOF. Let $\#_r x$, $\#_r y$ and $\#_r z$ be the numbers of $x_i$'s, $y_k$'s and $z_l$'s, respectively, that make up the configuration $\mathcal{R}$. By assumption, $\#_r x + \#_r y + \#_r z = r$ and, hence,

$$\#_r z = r - \#_r x - \#_r y \geq r - p - 2 = g - 2.$$

The claim being trivial if $\#_r z > g$, we consider the three remaining cases $\#_r z = g$, $g-1$ and $g-2$ separately. Now, $\#_r z = g$ implies $\#_r x + \#_r y = r - g = p \geq 2$; therefore, if there are no two $z_l$'s in one class, then one class must contain some $z_l$ together with an $x_i$ or a $y_k$. If $\#_r z = g-1$, then $\#_r x + \#_r y = r - g + 1 = p + 1$; since $p \geq 2$, at least one $x_i$ and one $y_k$ must belong to the configuration. A simple counting argument shows the claim in this case. Finally, if $\#_r z = g-2$, then $\#_r x + \#_r y = r - g + 2 = p + 2$, that is, all $x_i$'s and all $y_k$'s belong to the configuration. If all $z_l$'s form one-point classes, then the $x_i$'s and $y_k$'s must share the remaining two classes. If they are separated, then $\mathcal{R} = \mathcal{R}^\star$. In the opposite case, some class contains both an $x_i$ and a $y_k$. $\square$

For $0 \neq \mathbf{y} \in \mathbb{R}^d$, the rank-one matrix $\mathbf{y}\mathbf{y}^T$ has the simple eigenvalue $\|\mathbf{y}\|^2$. Therefore, $\det(\mathbf{I}_d + \mathbf{y}\mathbf{y}^T) = 1 + \|\mathbf{y}\|^2$, $\mathbf{y} \in \mathbb{R}^d$. Hence, if $\mathbf{A} \in \mathbb{R}^{d \times d}$ is positive definite, then we have

$$(16) \quad \det(\mathbf{A} + \mathbf{y}\mathbf{y}^T) = \det \sqrt{\mathbf{A}}(\mathbf{I}_d + \mathbf{A}^{-1/2}\mathbf{y}\mathbf{y}^T\mathbf{A}^{-1/2})\sqrt{\mathbf{A}} = (1 + \mathbf{y}^T\mathbf{A}^{-1}\mathbf{y}) \det \mathbf{A},$$

an equality that we will repeatedly use in the sequel.

In the following theorem we assume that, in the case of ties, an optimal solution is chosen for which the trace of the between-groups SSP matrix is minimum. This applies, in particular, to one-point clusters since these can be exchanged with discarded observations without any change of cost.

THEOREM 4.1 (Universal breakdown point of the TDC for the means).

(a) *If $n \geq r+1$ and $r \geq gd+2$, then the means remain bounded by a constant that depends only on the data as one observation is arbitrarily replaced.*

(b) *If $g \geq 2$ and $r \geq g+2$ (besides the standard assumption $r \geq gd+1$), then there is a data set such that one mean breaks down if two particular observations are suitably replaced.*

(c) *If $g \geq 2$, $n \geq r+1$, and $r \geq gd+2$, then $\beta_{\mathrm{mean}}(n,r,g) = \frac{2}{n}$.*

PROOF. (a) It suffices to show that an optimal configuration $\mathcal{R}^*$ discards a remote replacement. The mean vectors of all clusters of $\mathcal{R}^*$ will then remain within the convex hull of the data $\mathbf{x}_1^n$. Arguing by contradiction, let us assume that there is an optimal configuration $\mathcal{R}^*$ which contains the replacement $\mathbf{y}$ in a cluster $C \in \mathcal{R}^*$. If the replacement is far away, then this cluster must contain at least one other point since, otherwise, it would be exchanged with a discarded original point by the convention agreed upon just before the theorem. This point must, of course, be an original observation. It follows that $\mathbf{u} := \mathbf{y} - \mathbf{m}_C \to \infty$ as $\mathbf{y} \to \infty$. Since $r \geq gd+2$, any subset of $r$ elements, in particular the union of the members of $\mathcal{R}^*$, contains at least $r - 1 \geq gd + 1$ original points. Therefore, one cluster, $C_1$, contains a subset $E$ consisting of $d+1$ original points. The affine span of $E$ is the whole space and, whether $C_1 = C$ or not, we have

$$\mathbf{W}_{\mathcal{R}^*} \geq W_E(\mathbf{m}_{C_1}) + (\mathbf{y}_i - \mathbf{m}_C)(\mathbf{y}_i - \mathbf{m}_C)^T \geq W_E + \mathbf{u}\mathbf{u}^T.$$

Therefore, by (16),

$$\det \mathbf{W}_{\mathcal{R}^*} \geq \det(W_E + \mathbf{u}\mathbf{u}^T) = (1 + \|W_E^{-1/2}\mathbf{u}\|^2)\det W_E \xrightarrow[\mathbf{u}\to\infty]{} \infty.$$

This contradicts the fact that the maximum cost of any configuration that discards the replacement is finite.

(b) Let us construct a data set $D = \{\mathbf{x}_1, \ldots, \mathbf{x}_{r-g}, \mathbf{w}_1, \mathbf{w}_2, \mathbf{z}_1, \ldots, \mathbf{z}_{n-r+g-2}\}$. First note that $r - g \geq d+1$ (distinguish between $d=1$ and $d \geq 2$). Hence, we may choose $r-g$ elements $\{\mathbf{x}_1, \ldots, \mathbf{x}_{r-g}\} = F$ in general position with mean zero and SSP matrix $\mathbf{I}_d$.

We next use induction to construct the points $\mathbf{z}_1, \ldots, \mathbf{z}_{n-r+g-2}$ if $n - r + g > 2$. Suppose that $\mathbf{z}_1, \ldots, \mathbf{z}_l$ have already been constructed for some $l, 0 \leq l < n - r + g - 2$. Let $\mathcal{F}_l$ be the system of all hyperplanes $H$ spanned by $d$ of the points in $M_l = F \cup \{\mathbf{z}_1, \ldots, \mathbf{z}_l\}$ each. Since the set of all directions in $\mathbb{R}^d$ parallel to some $H \in \mathcal{F}_l$ is a $(d-1)$-dimensional subspace of $\mathbb{R}^d$, there is a direction not parallel to any of the hyperplanes. By running in such a direction, we find $\mathbf{z}_{l+1}$ so far from each hyperplane $H$ that

$$\text{(17)} \qquad \text{volume } P(\mathbf{z}_{l+1}, \mathbf{u}_1, \ldots, \mathbf{u}_d) \geq \sqrt{2(d+1)}$$

for all $\{\mathbf{u}_1, \ldots, \mathbf{u}_d\} \in \binom{M_l}{d}$ and such that

$$\text{(18)} \qquad (1 + \tfrac{1}{2}(\mathbf{z}_{l+1} - \mathbf{u})^T W_E^{-1}(\mathbf{z}_{l+1} - \mathbf{u})) \det W_E \geq 2$$

for all $\mathbf{u} \in M_l$ and all $E \in \binom{F \setminus \{\mathbf{u}\}}{d+1}$. After all $\mathbf{x}_i$'s and $\mathbf{z}_l$'s have been constructed, the two points $\mathbf{w}_1, \mathbf{w}_2$ are chosen in the centered unit ball so that $D$ is in general position. Irrespective of the optimal configuration, the $g$ estimated means are within the convex hull of $D$.

Now, mimicking the construction of one $\mathbf{z}_l$, we replace the two points $\mathbf{w}_1$ and $\mathbf{w}_2$ with a twin pair $\mathbf{y}_1 \neq \mathbf{y}_2$ such that $\|\mathbf{y}_2 - \mathbf{y}_1\| \leq 1$,

$$\text{(19)} \qquad \text{volume } P(E) \geq \sqrt{2(d+1)}$$

for all sets $E \in \binom{(D \setminus \{\mathbf{w}_1, \mathbf{w}_2\}) \cup \{\mathbf{y}_1, \mathbf{y}_2\}}{d+1}$ that contain at least one $\mathbf{y}_k$ except for $E = \{\mathbf{y}_1, \mathbf{y}_2\}$ if $d = 1$, and such that

$$\text{(20)} \qquad \det W_E (1 + \tfrac{1}{2}(\mathbf{y}_k - \mathbf{u})^T W_E^{-1}(\mathbf{y}_k - \mathbf{u})) \geq 2, \qquad k = 1, 2,$$

for all $\mathbf{u} \in D \setminus \{\mathbf{y}_1, \mathbf{y}_2\}$ and all $E \in \binom{F \setminus \{\mathbf{u}\}}{d+1}$.

We claim that the optimal configuration is $\mathcal{R}^\star = \{F, \{\mathbf{y}_1, \mathbf{y}_2\}, \{\mathbf{z}_1\}, \ldots, \{\mathbf{z}_{g-2}\}\}$. Indeed, by (16), its cost is

$$\begin{aligned}
\det \mathbf{W}_{\mathcal{R}^\star} &= \det(W_F + \tfrac{1}{2}(\mathbf{y}_2 - \mathbf{y}_1)(\mathbf{y}_2 - \mathbf{y}_1)^T) \\
&= \det(\mathbf{I}_d + \tfrac{1}{2}(\mathbf{y}_2 - \mathbf{y}_1)(\mathbf{y}_2 - \mathbf{y}_1)^T) \\
&= 1 + \tfrac{1}{2}\|\mathbf{y}_2 - \mathbf{y}_1\|^2 \leq \tfrac{3}{2}.
\end{aligned}$$

Moreover, Lemma 4.2 tells us that any clustering $\mathcal{R}$ not equivalent with $\mathcal{R}^\star$ (equivalent in the sense that some $\mathbf{z}_l$'s are exchanged) possesses some cluster $C$ (choose it of maximum size) containing either some pair $\{\mathbf{x}_i, \mathbf{y}_k\}$ or some $\mathbf{z}_l$ together with any other element. Let us denote these two elements by $\mathbf{a}$ and $\mathbf{b}$. If $C$ is of size at least $d+1$, then we choose a $(d+1)$-element subset $E \subseteq C$ containing $\{\mathbf{a}, \mathbf{b}\}$ and estimate

$$\det \mathbf{W}_\mathcal{R} \geq \det W_C \geq \det W_E \geq \frac{1}{d+1} \text{volume}^2 P(E) \geq 2$$

according to Lemma 4.1, (17) and (19).

Otherwise, we have $d \geq \#C \geq 2$ and there exists a cluster $R$ of size $\geq d+1$ which contains no pair $\{\mathbf{x}_i, \mathbf{y}_k\}$ and no $\mathbf{z}_l$. We have $R \neq C$ and from $d \geq 2$, it follows that $R \subseteq F$. Choosing a $(d+1)$-element subset $E \subseteq R$, we use (18), (20) and (16) to estimate

$$\det \mathbf{W}_{\mathcal{R}} \geq \det(W_E + W_{\{\mathbf{a},\mathbf{b}\}}) = \det W_E \det(1 + \tfrac{1}{2}(\mathbf{b}-\mathbf{a})^T W_E^{-1}(\mathbf{b}-\mathbf{a})) \geq 2.$$

In order to conclude the proof of part (b), it is sufficient to observe that the norm of the mean vector of the cluster $\{\mathbf{y}_1, \mathbf{y}_2\}$ may be arbitrarily large.

Part (c) follows immediately from (a) and (b). $\square$

4.2. *Restricted breakdown point and the separation property.* Theorem 4.1 says that the asymptotic universal breakdown value for the means is zero; this is a negative result and somewhat unsatisfactory in the framework of a trimming algorithm. One reason is the strength of the universal breakdown value. We may rescue the situation by introducing a relaxed version of it, the *restricted breakdown value* $\beta(\delta, \mathcal{K})$ w.r.t. a subclass $\mathcal{K} \subseteq \mathcal{A}$ of data sets admissible for $\delta$. It lies between the individual and the universal breakdown values, Section 4.1, and we define it as the minimum fraction of replacements that cause $\delta$ to break down with some $\mathbf{x}_1^n \in \mathcal{K}$,

$$\beta(\delta, \mathcal{K}) = \min_{\mathbf{x}_1^n \in \mathcal{K}} \beta(\delta, \mathbf{x}_1^n).$$

The universal breakdown value is just $\beta(\delta, \mathcal{A})$. The restricted breakdown value, too, is a characteristic of an estimator. It provides information on the structure that a data set must have so that the estimator still acts reasonably after contamination.

Let us now compute the restricted breakdown value of (TDC) for the mean values w.r.t. a certain class of data sets that we describe next. It is necessary to first introduce some notation. Let $\mathcal{P} = \{P_1, \ldots, P_g\}$ be a partition of some data set $S$, let $\mathcal{E}_j^{\mathcal{P}}$ be the set of all mean values of nonempty subsets of $P_j$, $1 \leq j \leq g$, let $\mathcal{S}^{\mathcal{P}}$ be the set of all subconfigurations of $\mathcal{P}$ comprising at most $r$ elements and possessing at least one cluster of size $\geq d+1$, and let $\mathcal{W}^{\mathcal{P}}$ be the set of all pooled SSP matrices generated by elements of $\mathcal{S}^{\mathcal{P}}$. Given $g \geq 2$ and $u \geq 1$, we define $k_{g,u} = \lceil \frac{\max\{2r-n, (g-1)gd+1, n-u+1\}}{(g-1)g} \rceil$ $(> d)$ and $\mathcal{K}_{g,u}$ as the system of all $d$-dimensional data sets $S$ of length $n$ in general position that have the *separation property*

$S$ possesses a partition $\mathcal{P}$ in $g$ subsets of sizes at least $u$ $(\leq \lfloor n/g \rfloor)$ such that

$$(21) \min_{\substack{\mathbf{W} \in \mathcal{W}^{\mathcal{P}} \\ \mathbf{m}_j \in \mathcal{E}_j^{\mathcal{P}} \\ \mathbf{m}_k \in \mathcal{E}_k^{\mathcal{P}} \\ j \neq k}} (\mathbf{m}_k - \mathbf{m}_j)^T \mathbf{W}^{-1} (\mathbf{m}_k - \mathbf{m}_j) > 2 \cdot \frac{\max_{\mathbf{W} \in \mathcal{W}^{\mathcal{P}}} \det \mathbf{W}}{\min_{\substack{C \in \binom{P_j}{k_{g,u}} \\ 1 \leq j \leq g \\ \#P_j \geq k_{g,u}}} \det W_C}.$$

Both sides of this estimate are invariant with respect to location and scale and their quotient describes a measure of *validity* of the partition $\mathcal{P}$, a combination of cluster *compactness* ($\det \mathbf{W}$) and cluster *separation* (the Mahalanobis distance squared). A great number of such indices are widely in use for assessing the quality of a partition; see Bezdek, Keller, Krisnapuram and Pal (1999).

A data set satisfying condition (21) possesses a marked cluster structure. Note that the left-hand side of (21) increases as the different clusters in $\mathcal{P}$ are moved away from each other, whereas the right-hand side remains unchanged under this operation. It is, therefore, easy to construct examples of data sets that satisfy the separation property.

Note that the classes $\mathcal{K}_{g,u}$ decrease as $u$ increases. The SSP matrix $\mathbf{W}_\mathcal{P}$ is larger than all matrices in $\mathcal{W}^\mathcal{P}$ w.r.t. the positive semi-definite ordering; therefore, substituting $\min_{\mathbf{m}_j,\mathbf{m}_k}(\mathbf{m}_k - \mathbf{m}_j)^T \mathbf{W}_\mathcal{P}^{-1}(\mathbf{m}_k - \mathbf{m}_j)$ for the left-hand side of (21) and $\det \mathbf{W}_\mathcal{P}$ for $\max_{\mathbf{W} \in \mathcal{W}^\mathcal{P}} \det \mathbf{W}$ defines a narrower class. It is easier to verify this condition than (21).

LEMMA 4.3. *Let $S \in \mathcal{K}_{g,1}$ (with associated partition $\mathcal{P}$) and let $\mathcal{R}$ be a partition of some finite subset $R \subseteq \mathbb{R}^d$ such that:*

(i) *some cluster $R_k \in \mathcal{R}$ contains elements of at least two different $P_j$'s;*
(ii) *there are a cluster $R_l \in \mathcal{R}$ and some $j \in 1..g$ such that $\#(R_l \cap P_j) \geq k_{g,u}$.*

*Then we have $\det \mathbf{W}_\mathcal{R} > \max_{\mathbf{W} \in \mathcal{W}^\mathcal{P}} \det \mathbf{W}$.*

PROOF. Without loss of generality, the two $P_j$'s appearing in (i) are $P_1$ and $P_2$. We first consider the case $k = l$. Putting $a_j = \#(R_k \cap P_j)$, we use Lemma A.3 to estimate

$$\mathbf{W}_\mathcal{R} \geq W_{R_k} = \sum_{x \in R_k} (x - \mathbf{m}_{R_k})(x - \mathbf{m}_{R_k})^T$$

$$= \sum_{j=1}^{g} \sum_{x \in R_k \cap P_j} (x - \mathbf{m}_{R_k})(x - \mathbf{m}_{R_k})^T$$

$$= \sum_{j=1}^{g} W_{R_k \cap P_j} + \sum_{1 \leq j < h \leq g} \frac{a_j a_h}{\#R_k}(\mathbf{m}_{R_k \cap P_j} - \mathbf{m}_{R_k \cap P_h})(\mathbf{m}_{R_k \cap P_j} - \mathbf{m}_{R_k \cap P_h})^T$$

$$= \mathbf{W}_{R_k \cap \mathcal{P}} + \sum_{1 \leq j < h \leq g} \frac{a_j a_h}{\#R_k}(\mathbf{m}_{R_k \cap P_j} - \mathbf{m}_{R_k \cap P_h})(\mathbf{m}_{R_k \cap P_j} - \mathbf{m}_{R_k \cap P_h})^T;$$

here we have abbreviated $R_k \cap \mathcal{P} = (R_k \cap P_1, \ldots, R_k \cap P_g)$. Applying Lemma A.1(b) with $A = \mathbf{W}_{R_k \cap \mathcal{P}}$ and $y_{jh} := \sqrt{\frac{a_j a_h}{\#R_k}}(\mathbf{m}_{R_k \cap P_j} - \mathbf{m}_{R_k \cap P_h})$, $1 \leq j < h \leq g$, we

infer

$$\det \mathbf{W}_{\mathcal{R}} \geq \left(1 + \sum_{1 \leq j < h \leq g} \frac{a_j a_h}{\#R_k}(\mathbf{m}_{R_k \cap P_j} - \mathbf{m}_{R_k \cap P_h})^T \mathbf{W}_{R_k \cap \mathcal{P}}^{-1}(\mathbf{m}_{R_k \cap P_j} - \mathbf{m}_{R_k \cap P_h})\right)$$
$$\times \det \mathbf{W}_{R_k \cap \mathcal{P}}$$
$$> \left(\sum_{1 \leq j < h \leq g} \frac{a_j a_h}{\#R_k}\right) \min_{\substack{R_k \cap P_j, R_k \cap P_h \neq \varnothing \\ 1 \leq j < h \leq g}} (\mathbf{m}_{R_k \cap P_j} - \mathbf{m}_{R_k \cap P_h})^T$$
$$\times \mathbf{W}_{R_k \cap \mathcal{P}}^{-1}(\mathbf{m}_{R_k \cap P_j} - \mathbf{m}_{R_k \cap P_h}) \det \mathbf{W}_{R_k \cap \mathcal{P}}$$
$$\geq \frac{1}{2} \min_{\substack{R_k \cap P_j, R_k \cap P_h \neq \varnothing \\ 1 \leq j < h \leq g}} (\mathbf{m}_{R_k \cap P_j} - \mathbf{m}_{R_k \cap P_h})^T \mathbf{W}_{R_k \cap \mathcal{P}}^{-1}(\mathbf{m}_{R_k \cap P_j} - \mathbf{m}_{R_k \cap P_h})$$
$$\times \det \mathbf{W}_{R_k \cap \mathcal{P}};$$

the last inequality follows from Lemma A.4 and (i). Since $R_k \cap \mathcal{P} \in \mathcal{S}^{\mathcal{P}}$ by (ii) and since $\mathbf{m}_{R_k \cap P_h} \in \mathcal{E}_h^{\mathcal{P}}$, $1 \leq h \leq 2$, we may apply the lower bound (21) to the last line above to obtain

$$\det \mathbf{W}_{\mathcal{R}} > \frac{\max_{\mathbf{W} \in \mathcal{W}^{\mathcal{P}}} \det \mathbf{W}}{\min_{\substack{C \in \binom{P_j}{k_{g,u}} \\ 1 \leq j \leq g \\ \#P_j \geq k_{g,u}}} \det W_C} \det \mathbf{W}_{R_k \cap \mathcal{P}} \geq \max_{\mathbf{W} \in \mathcal{W}^{\mathcal{P}}} \det \mathbf{W}.$$

If $k \neq l$, we start with

$$\mathbf{W}_{\mathcal{R}} \geq W_{R_l} + W_{R_k}$$
$$\geq \mathbf{W}_{R_l \cap \mathcal{P}} + \sum_{1 \leq j < h \leq g} \frac{a_j a_h}{\#R_k}(\mathbf{m}_{R_k \cap P_j} - \mathbf{m}_{R_k \cap P_h})(\mathbf{m}_{R_k \cap P_j} - \mathbf{m}_{R_k \cap P_h})^T.$$

The remainder of the proof is similar, as before. □

If a data set meets the separation property, then (TDC) is much more robust than predicted by Theorem 4.1.

THEOREM 4.2 (Restricted breakdown point of the TDC for the means). *Let $g \geq 2$, let $r > (g-1)gd$, and let $u > n - r$ be an integer.*

(a) *The restricted breakdown value of (TDC) for the mean values w.r.t. $\mathcal{K}_{g,u}$ satisfies $\beta_{\mathrm{mean}}(n, r, g, \mathcal{K}_{g,u}) \geq \frac{1}{n} \min\{n - r + 1, r - (g-1)gd, r + u - n\}$.*
(b) *If*

(i) $2r - n > (g-1)gd$ *and*
(ii) $u > 2(n - r)$,

*then $\beta_{\mathrm{mean}}(n, r, g, \mathcal{K}_{g,u}) = \frac{1}{n}(n - r + 1)$.*

PROOF. (a) Let $S \in \mathcal{K}_{g,u}$ with partition $\mathcal{P}$ and let $M$ be any set obtained from $S$ by modifying at most $\rho = \min\{n - r, r - (g-1)gd - 1, r + u - n - 1\}$ elements. Our proof proceeds in several steps.

($\alpha$) Any configuration $\mathcal{R}$ in $M$ has at least one cluster with $d+1$ *original* observations:

Indeed, by definition of $\rho$, $R = \bigcup \mathcal{R}$ has at least $r - \rho > (g-1)gd \geq gd$ original observations and the claim follows from the pigeon hole principle.

Let $\mathcal{R}$ now be the optimal configuration of $M$. We will show that the mean values of all clusters of $\mathcal{R}$ are bounded by a number that depends solely on the original data $S$.

($\beta$) $\det \mathbf{W}_\mathcal{R}$ is bounded by a number that depends only on $S$:

In fact, we have

$$\det \mathbf{W}_\mathcal{R} \leq \max_{\mathbf{W} \in \mathcal{W}^\mathcal{P}} \det \mathbf{W}. \tag{22}$$

Indeed, let $R' \subseteq M$ consist of $r$ original points and let $\mathcal{R}' = R' \cap \mathcal{P} = (R' \cap P_1, \ldots, R' \cap P_g)$; then $\mathcal{R}' \in \mathcal{S}^\mathcal{P}$, $\mathbf{W}_{\mathcal{R}'} \in \mathcal{W}^\mathcal{P}$ and $\det \mathbf{W}_\mathcal{R} \leq \det \mathbf{W}_{\mathcal{R}'}$ by optimality.

($\gamma$) If $R_j$ contains $d+1$ or more original observations, then its mean $m_j$ is bounded by a number that depends only on $S$:

To this end, define

$$\lambda_{\max}(S) := \max\{\lambda | \lambda \text{ eigenvalue of } W_C, C \subseteq S \text{ and } \#C \geq d+1\},$$
$$b_{\min}(S) := \min\{\det W_C | C \subseteq S, \#C \geq d+1\}.$$

These quantities are bounded above and below (away from 0) and depend only on $S$. Now, by Steiner's identity,

$$\mathbf{W}_\mathcal{R} \geq \sum_{\mathbf{x} \in R_j \cap S} (\mathbf{x} - \mathbf{m}_{R_j})(\mathbf{x} - \mathbf{m}_{R_j})^T$$
$$= W_{R_j \cap S} + \#(R_j \cap S)(\mathbf{m}_{R_j} - \mathbf{m}_{R_j \cap S})(\mathbf{m}_{R_j} - \mathbf{m}_{R_j \cap S})^T$$
$$\geq W_{R_j \cap S} + (\mathbf{m}_{R_j} - \mathbf{m}_{R_j \cap S})(\mathbf{m}_{R_j} - \mathbf{m}_{R_j \cap S})^T.$$

Hence, by (16) and the assumption made on $R_j$,

$$\det \mathbf{W}_\mathcal{R} \geq \det W_{R_j \cap S}(1 + (\mathbf{m}_{R_j} - \mathbf{m}_{R_j \cap S})^T W_{R_j \cap S}^{-1}(\mathbf{m}_{R_j} - \mathbf{m}_{R_j \cap S}))$$
$$> \frac{b_{\min}(S)}{\lambda_{\max}(S)} \|\mathbf{m}_{R_j} - \mathbf{m}_{R_j \cap S}\|^2$$

and the claim now follows from ($\beta$).

($\delta$) If $R_j$ contains between one and $d$ original observations, then its mean $m_j$ is bounded by a number that depends only on $S$:

By ($\alpha$), there is a cluster $k \neq j$ containing $d+1$ original observations. We have

$$\mathbf{W}_{\mathcal{R}} \geq W_{R_k} + W_{R_j}$$
$$\geq W_{R_k \cap S} + W_{R_j \cap S} + \#(R_j \cap S)(\mathbf{m}_{R_j} - \mathbf{m}_{R_j \cap S})(\mathbf{m}_{R_j} - \mathbf{m}_{R_j \cap S})^T$$
$$\geq W_{R_k \cap S} + (\mathbf{m}_{R_j} - \mathbf{m}_{R_j \cap S})(\mathbf{m}_{R_j} - \mathbf{m}_{R_j \cap S})^T$$

by assumption on $R_j$; hence,

$$\det \mathbf{W}_{\mathcal{R}} \geq \det W_{R_k \cap S}(1 + (\mathbf{m}_{R_j} - \mathbf{m}_{R_j \cap S})^T W_{R_k \cap S}^{-1}(\mathbf{m}_{R_j} - \mathbf{m}_{R_j \cap S}))$$

and the proof terminates as that of ($\gamma$).

In view of ($\gamma$) and ($\delta$), the proof of part (a) will be finished if we show that

($\varepsilon$) each $R_j$ contains at least one original point.

Assume, on the contrary, that $\mathcal{R}$ contains some cluster that consists solely of replacements. We show that, as a consequence, $R$ and $\mathcal{R}$ must satisfy the hypotheses of Lemma 4.3. By definition of $\rho$, $R$ has at least $r - (r + u - n - 1) = n + 1 - u$ original elements; that is, $\sum_{j=1}^{g} \#R \cap P_j > n - u$. Taking into account that each $P_j$ has at least $u$ elements, a simple counting argument shows that each of the $g$ sets $P_j$ intersects $R$. By assumption, there is some $R_h$ omitted by all $P_j$'s and the pigeon hole principle shows Lemma 4.3(i). Again by assumption, the original observations in $R$ are distributed over at most $(g-1)g$ (disjoint) sets of the form $R_l \cap P_j$; the definition of $k_{g,u}$ and another application of the pigeon hole principle show that some set $R_l \cap P_j$ contains at least $k_{g,u}$ original observations; this is Lemma 4.3(ii).

The conclusion of Lemma 4.3 contradicts (22), which completes the proof of part (a).

(b) The individual breakdown value of (TDC) for the mean vectors w.r.t. any admissible data set $\mathbf{x}_1^n$ is $\leq \frac{n-r+1}{n}$. Indeed, let $M$ be a set obtained from $\mathbf{x}_1^n$ by modifying at least $n - r + 1$ of its elements. Then each subset of $M$ of size $r$ contains at least $r - (n - (n - r + 1)) = 1$ replacements. Part (b) now follows from (a), (i) and (ii). $\square$

In the case $g = 1$, (TDC) reduces to Rousseeuw's MCD; see Corollary 2.2. Rousseeuw (1985) proves that the asymptotic breakdown value of MCD with $r = \lceil (1-\alpha)n \rceil$ is $\alpha$ if $\alpha < 0.5$; see also Lopuhaä and Rousseeuw (1991), where the analogous estimator MVE is treated in more detail. This is in harmony with the foregoing theorem, although it is not a corollary of it: the separation property and part ($\varepsilon$) of its proof are not applicable if $g = 1$.

Theorem 4.2(b) asserts that the algorithm can withstand exactly the number of outliers generated by the model 2.1 if the hypotheses of (b) are satisfied.

EXAMPLE 4.1. By way of discussing Theorems 4.1, 4.2 and the separation property, it is interesting to take a look at an instructive example. Let us consider the one-dimensional data set $x_1,\ldots,x_{10}$ shown in Figure 1 with gap $a>1$. Its "natural" partition has the $g=2$ clusters $R_1=\{x_1,\ldots,x_5\}$ and $R_2=\{x_6,\ldots,x_{10}\}$, and the means are 0 and $a+4$. It is reasonable to choose this partition for $\mathcal{P}$ and $u=5$. The assumptions of both theorems are met with the parameter $r=8$ (i.e., the algorithm discards two observations). The first theorem says that (TDC) will resist one arbitrary outlier for all $a>1$, whereas the second promises that it will tolerate even two such outliers if $a>20$. There is actually a transition from the breakdown value 0.2 to 0.3 at a much smaller value of $a$. Let us compute this critical value. A decisive pair of replacements is $(x_7, x_8)$. As we replace these two observations by very large, close numbers $x_7'$, $x_8'$ with SSP $\varepsilon$, two configurations compete for optimality: the configurations $R_1'=\{x_1,\ldots,x_6\}$, $R_2'=\{x_7',x_8'\}$ ($x_9$ and $x_{10}$ removed, one mean breaks down) and $R_1''=\{x_1,\ldots,x_5\}$, $R_2''=\{x_6,x_9,x_{10}\}$ (replacements removed, means are 0 and $a+\frac{13}{3}$). The first has SSP $10+\frac{5}{6}(a+2)^2+\varepsilon$, which is smaller than that of the second, $\frac{56}{3}$, if $a<\sqrt{\frac{52}{5}}-2\approx 1.225$. Hence, this is the critical parameter that separates the breakdown values 0.2 and 0.3. This indicates that (TDC) is actually more robust than predicted by Theorem 4.2, let alone Theorem 4.1.

We next show that, if $n>r$ and if $r/n$ is large enough, the TDC is robust w.r.t. the SSP matrix. We actually show that it breaks down simultaneously for each data set as the fraction of bad outliers slightly exceeds $1-r/n$.

THEOREM 4.3 (Universal breakdown point of the TDC for the pooled SSP matrix).

(a) *Suppose that $2r \geq n + g(d+1)$. If at most $n-r+g-1$ points of the data set $D$ are replaced in an arbitrary way, then the eigenvalues of the SSP*

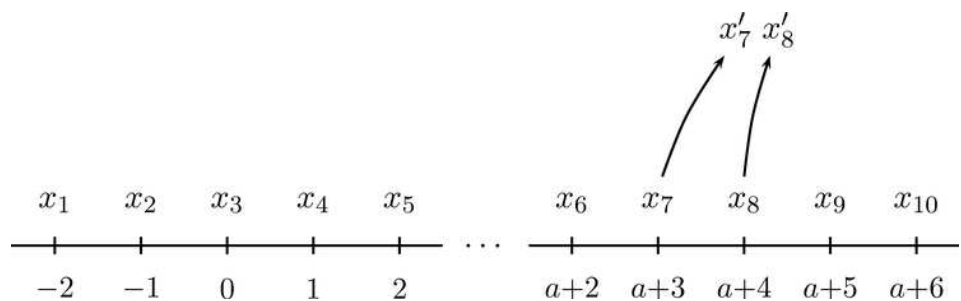

FIG. 1. 1D data set, replacements $x_7'$, $x_8'$; breakdown at 2 replacements if $a<1.225$ and at 3 if $a \geq 1.225$.

*matrix of any admissible configuration remain bounded away from zero by a constant that depends only on $D$ and $d$.*

(b) *Suppose that $2r \geq n + g(d+1)$. If at most $n - r + g - 1$ points of $D$ are replaced in an arbitrary way, then the eigenvalues of the SSP matrix of the optimal configuration remain bounded by a constant that depends only on $D$ and $d$.*

(c) *Given $t > 0$, $n - r + g$ elements of any $D$ may be replaced so that the largest eigenvalue of the SSP matrix exceeds $t$.*

(d) *If $2r \geq n + g(d+1)$, then $\beta_{\mathrm{SSP}}(n,r,g) = \frac{n-r+g}{n}$.*

PROOF. We begin with two remarks. (i) Let $D$ be any subset of $\mathbb{R}^d$ of size $n$ in general position and let $W_E$ be the SSP matrix of a subset $E \subseteq D$. Since $D$ is in general position, the number

$$\alpha = \min_{E \in \binom{D}{d+1}} \lambda_{\min}(W_E)$$

is strictly positive and depends on $D$ and $d$, only.

(ii) Let $M$ be any set obtained from $D$ by modifying at most $n - r + g - 1$ elements and let $R$ be any subset of $M$ consisting of $r$ elements. If $2r \geq n + g(d+1)$, then $R$ contains at least $r - (n - r - g - 1) = 2r - n - g + 1 \geq gd + 1$ original observations in $D$. By the pigeon hole principle, any clustering of $R$ in $g$ parts has at least one member $C$ that contains $d + 1$ such points. We now prove (a)–(d).

(a) The least eigenvalue of the SSP matrix of $C$ in (ii) is $\geq \alpha$ and, by Section 4.1 the same is true for the SSP matrix of any admissible clustering of the modified set $M$.

(b) Since $r - g + 1$ points in $D$ remain unchanged by the replacement, the modified set $M$ has an admissible configuration $\mathcal{M}$ with one cluster consisting of $r - g + 1$ original data and $g - 1$ one-point clusters. Hence, the SSP matrix of the optimal configuration $\mathcal{M}^\star$ of $M$ cannot have a determinant larger than that of $\mathcal{M}$, that is,

$$\det \mathbf{W}_{\mathcal{M}^\star} \leq \det \mathbf{W}_{\mathcal{M}}.$$

The SSP matrix $\mathbf{W}_{\mathcal{M}}$ is that of the $r - g + 1$ original points and, hence, depends only on $D$. Therefore, its determinant is bounded by a constant $\gamma$ that again depends only on $D$. By (ii), at least one cluster of $\mathcal{M}^\star$ contains at least $d + 1$ original points. The claim, therefore, follows from the estimates

$$\lambda_{\max}(\mathbf{W}_{\mathcal{M}^\star})\alpha^{d-1} \leq \lambda_{\max}(\mathbf{W}_{\mathcal{M}^\star})\lambda_{\min}(\mathbf{W}_{\mathcal{M}^\star})^{d-1} \leq \det \mathbf{W}_{\mathcal{M}^\star} \leq \det \mathbf{W}_{\mathcal{M}} \leq \gamma.$$

(c) Modify $D$ by $\geq n - r + g$ replacements that are at a distance $\geq 2t$ from each original observation and from each other to obtain a set $M$. Let $\mathcal{M}^\star$ be its optimal clustering. Clearly, any subset of $r$ elements of $M$ contains at least $g$ replacements. Moreover, there is a cluster $C$, $\#C \geq 2$, that contains at least one replacement. Indeed, if no cluster contains two replacements, then each contains exactly one and, since $r \geq gd+1$, one cluster contains at least two elements. Now, if $\mathbf{x}$ is a replacement and $\mathbf{y}$ another element in $C$, then the trace of the SSP matrix of $C$ is at least

$$\operatorname{tr} W_C \geq \operatorname{tr} \left\{ \left(\mathbf{x} - \frac{\mathbf{x}+\mathbf{y}}{2}\right)\left(\mathbf{x} - \frac{\mathbf{x}+\mathbf{y}}{2}\right)^T + \left(\mathbf{y} - \frac{\mathbf{x}+\mathbf{y}}{2}\right)\left(\mathbf{y} - \frac{\mathbf{x}+\mathbf{y}}{2}\right)^T \right\}$$
$$= \frac{\|\mathbf{x} - \mathbf{y}\|^2}{2} \geq t.$$

Therefore, $\operatorname{tr} \mathbf{W}_{\mathcal{M}^\star} \geq \operatorname{tr} W_C \geq t$ (see 4.1) and one eigenvalue must exceed $t/d$.

Claim (d) follows from (a)–(c). $\square$

It may be astonishing that the TDC should take $g-1$ replacements even if $r = n$. However, isolated replacements may be treated as one-point clusters in the optimal configuration; in this case, the number of clusters formed by the original data is reduced. Replacements that huddle together may form one cluster. In both cases, the SSP matrix is not completely destroyed.

**5. Simulation studies.** In order to assess the performance of the proposed algorithm, we have implemented it as a C++ program for various dimensions, sizes, numbers and positions of clusters, as well as numbers of outliers. The first simulation study illustrates how, by varying the input parameter $r$ (the assumed number of regular observations) of our algorithm, one can roughly control the amount of contaminations contained in the data set.

The symbol $\mathbf{e}_i \in \mathbb{R}^d$, $i \in 1..d$, stands for the $i$th unit vector. As usual, the symbol $\chi_d^2(\alpha)$ denotes the $\alpha$-quantile of the $\chi^2$-distribution with $d$ degrees of freedom. We consider, in dimension $d = 8$, the $2d$ normally distributed populations $N_d(\boldsymbol{\mu}_j, \mathbf{V})$, $j \in 1..2d$, with common covariance matrix $\mathbf{V}$, diagonal with entries $(1.0, 1.2, 1.4, 1.6, 1.8, 2.0, 2.2, 9.0)$ and means

$$(23) \quad \boldsymbol{\mu}_j = \begin{cases} -\sqrt{\dfrac{\mathbf{V}((j+1)/2, (j+1)/2)\chi_d^2(\alpha)}{2}} \mathbf{e}_{(j+1)/2}, & j \text{ odd}, \\ \sqrt{\dfrac{\mathbf{V}(j/2, j/2)\chi_d^2(\alpha)}{2}} \mathbf{e}_{j/2}, & j \text{ even}, \end{cases}$$

$\alpha \in \{0.95, 0.99, 0.999, 0.999999\}$. That is, the means of the various clusters lie on the $2d$ axial directions of $\mathbb{R}^d$; those of the $j$th and $(j+1)$st clusters, $j$ odd, lie on the same axis but in opposite directions viewed from the origin. The means (23) assure that the squared Mahalanobis distance of two cluster centers on the same axis is $2\chi_d^2(\alpha)$, whereas it is $\chi_d^2(\alpha)$ in the opposite case. The values $\alpha = 0.95$ and $0.99$ give rise to heavily and moderately overlapped clusters, while $\alpha = 0.999$ and $0.999999$ mean better and good separation, respectively. We generate 100 observations from each cluster. Thus, we obtain a total of $r = 200\,d$ (regular) observations. Additionally, in our first essay, we contaminate the data with $22\,d$ outliers arranged in shells around the cluster centers. The square of the Mahalanobis distance from each contamination to the closest cluster center is $\chi_d^2(\beta)$, $\beta \in \{0.999, 0.999999\}$; see Figure 2.

Since we consider four $\alpha$'s and two $\beta$'s, these specifications define eight different cases. To each of them we apply the algorithm described in Section 3 four times, namely, with the a priori numbers of regular observations $r = n, 0.95n, 0.9n$ and $0.85n$. More precisely, we apply the multistart method with up to 2000 random initial configurations based on the method 3.3(a) and iterate reduction steps until convergence is reached. The 32 rows in Table 1 show, for each choice of $\alpha$, $\beta$ and $r$, the fractions of estimated regular observations whose squared Mahalanobis distances to their estimated cluster centers (w.r.t. the estimated common covariance matrix) are larger than a given percentile $\chi_d^2(\gamma)$, $\gamma \in \{0.95, 0.975, 0.99, 0.999\}$. In the rows corresponding to the correct fraction of outliers ($\approx 10\%$), these are expected to match the correct tail probabilities $1 - \gamma$ shown on the top of Table 1. This heuristic (akin to a $\chi^2$ goodness-of-fit test) for estimating the number of outliers is suggested by a similar heuristic applied in robust discriminant analysis by Gather and Kale [(1988), Section 3] and Ritter and Gallegos (1997). The method slightly underestimates the number of outliers. We also compare the theoretical populations with those defined by the estimated clusters. The Bhattacharyya distance between two normal distributions is

$$d_{\text{Bhatt}}(N_d(\boldsymbol{\mu}_1, \mathbf{V}_1), N_d(\boldsymbol{\mu}_2, \mathbf{V}_2))$$
$$= 1 - \sqrt{\frac{\sqrt{\det \mathbf{V}_1 \det \mathbf{V}_2}}{\det((\mathbf{V}_1 + \mathbf{V}_2)/2)}} \exp\left(-\frac{1}{4}(\boldsymbol{\mu}_2 - \boldsymbol{\mu}_1)^T (\mathbf{V}_1 + \mathbf{V}_2)^{-1}(\boldsymbol{\mu}_2 - \boldsymbol{\mu}_1)\right),$$

a number in the unit interval. A measure for the quality of the estimates is the minimum over all matchings $\sigma \in \mathcal{S}_{2d}$ between estimated and real classes of the maximum Bhattacharyya distance over the $2d$ matched pairs:

$$\min_{\sigma \in \mathcal{S}_{2d}} \max_{j \in 1..2d} d_{\text{Bhatt}}\left(N_d\left(\mathbf{m}_{\sigma(j)}, \frac{1}{r}\mathbf{W}\right), N_d(\boldsymbol{\mu}_j, \mathbf{V})\right).$$

The results are shown in Table 2 for $d = 2, 4, 8$. For $d = 2$ and $\alpha = 0.95$, the original clusters are not recovered by the algorithm since the 400 regular data points are too homogeneous; there is no reasonable matching. We tested also scenarios with *diffuse* outliers generated from $N_d(\boldsymbol{\mu}, v \cdot \mathbf{I}_d)$, $d = 2, 4, 8$, $\boldsymbol{\mu} \in \mathbb{R}^d$ and $v \geq 16$. The case $\boldsymbol{\mu} = 0$, $v = 16$ is the most demanding since the variance is already quite close to that of the last coordinate. Even in this case, only about 10% of the rejected elements are (extreme) regular observations and the generated clusters are well rediscovered.

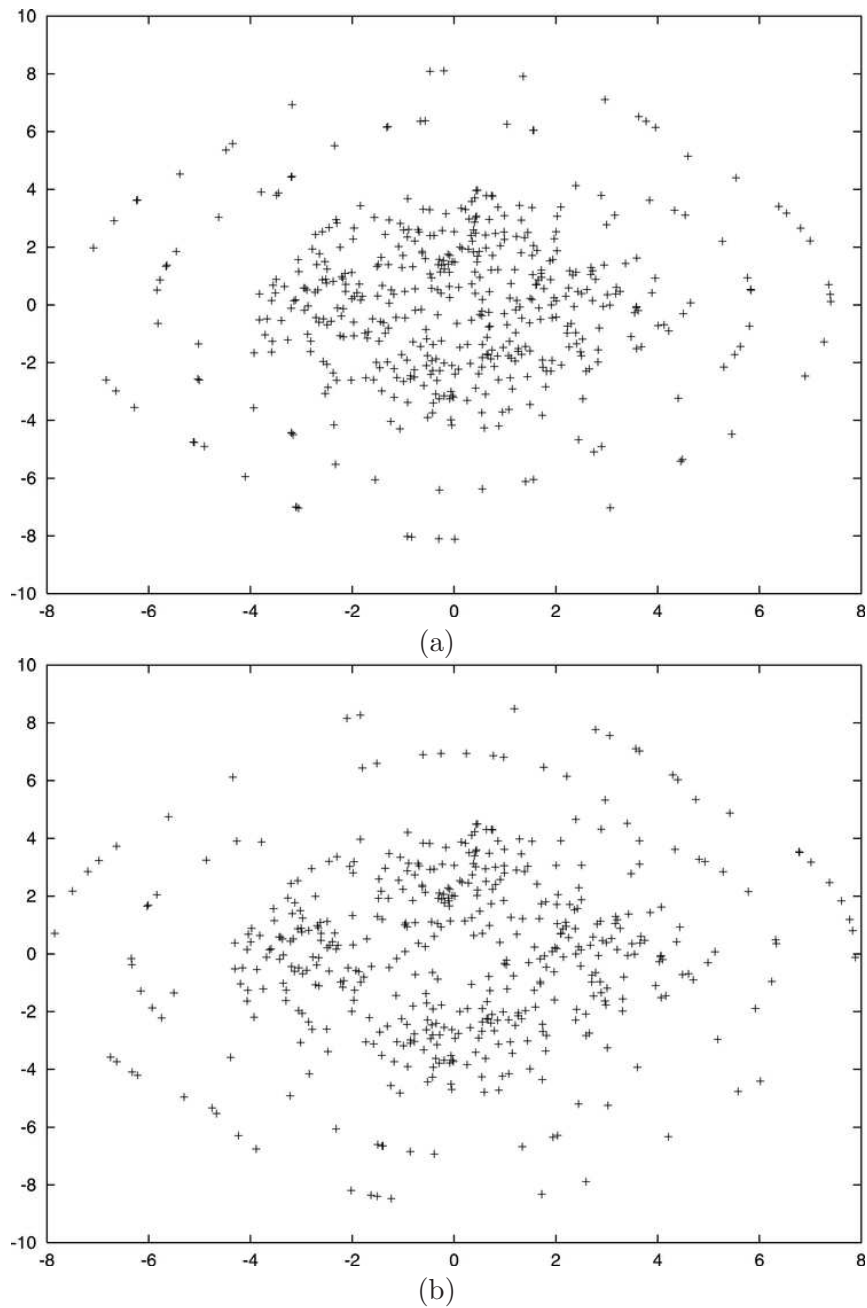

Fig. 2. *Visualization of four clusters in two dimensions. The Mahalanobis distance between each pair of the four cluster centers depends on the value $\alpha$ as specified by* (23). *The outliers are arranged in shells. The inner and outer shells correspond to $\beta = 0.999$ and $\beta = 0.999999$, respectively, where $\beta$ defines the ellipsoids of equal concentration on which the contaminations lie. See also the text.* (a) $\alpha = 0.99$, (b) $\alpha = 0.999$.

Finally, the number $L$ of reduction steps until convergence in one iteration is about 15 for $d=2$ and 22 for $d=4,8$ with standard deviations about 7. One reduction step takes no longer than 0.004, 0.01 and 0.06 seconds, respectively, on a 1 GHz processor. These figures are essentially independent of the trimming parameter $r$ of the algorithm.

We do not contend that the proposed algorithm responds to each clustering situation. In fact, the presented model is meant as a possible component for outlier handling in a comprehensive clustering strategy. One of the main purposes of the paper is a contribution to computing breakdown values. Nevertheless, in our experience, the algorithm works well in situations where the model assumptions are satisfied: clusters of about the same shape, scale and size, and outliers sufficiently scattered in space (not concentrated close to one or a few points).

## APPENDIX: TECHNICAL PRELIMINARIES

In this Appendix we prepare some tools for the proofs of the theorems. Some of the statements are of interest on their own.

LEMMA A.1. *Let $d \geq 2$ and let $\mathbf{A} \in \mathbb{R}^{d \times d}$ be positive definite.*

(a) *For any positive semi-definite matrix $\mathbf{C} \in \mathbb{R}^{d \times d}$, we have*

$$\det(\mathbf{A} + \mathbf{C}) \geq (1 + \operatorname{tr}(\mathbf{A}^{-1}\mathbf{C})) \det \mathbf{A} + \det \mathbf{C}.$$

(b) *If $y_1, \ldots, y_k \in \mathbb{R}^d$, then we have*

$$\det\left(\mathbf{A} + \sum_h y_h y_h^T\right) \geq \left(1 + \sum_h y_h^T A^{-1} y_h\right) \det \mathbf{A}.$$

PROOF. (a) From $\mathbf{A} + \mathbf{C} = \mathbf{A}^{1/2}(\mathbf{I}_d + \mathbf{A}^{-1/2}\mathbf{C}\mathbf{A}^{-1/2})\mathbf{A}^{1/2}$, we infer

$$\det(\mathbf{A} + \mathbf{C}) = \det(\mathbf{I}_d + \mathbf{A}^{-1/2}\mathbf{C}\mathbf{A}^{-1/2}) \det \mathbf{A}. \tag{24}$$

If $\lambda_1, \lambda_2, \ldots, \lambda_d$ are the eigenvalues of $\mathbf{A}^{-1/2}\mathbf{C}\mathbf{A}^{-1/2}$, then the eigenvalues of $\mathbf{I}_d + \mathbf{A}^{-1/2}\mathbf{C}\mathbf{A}^{-1/2}$ are $1 + \lambda_1, \ldots, 1 + \lambda_d$ and the claim follows from (24) and

$$\det(\mathbf{I}_d + \mathbf{A}^{-1/2}\mathbf{C}\mathbf{A}^{-1/2}) = \prod_{j=1}^d (1 + \lambda_j) \geq 1 + \sum_{j=1}^d \lambda_j + \prod_{j=1}^d \lambda_j$$

$$= 1 + \operatorname{tr}(\mathbf{A}^{-1/2}\mathbf{C}\mathbf{A}^{-1/2}) + \det \mathbf{A}^{-1/2}\mathbf{C}\mathbf{A}^{-1/2}.$$

Part (b) is an immediate consequence of (a). □

LEMMA A.2. $0 \leq \mathbf{A} \leq \mathbf{B}$ *implies* $\det \mathbf{A} \leq \det \mathbf{B}$. *If $\mathbf{B} > 0$, then there is equality if and only if $\mathbf{A} = \mathbf{B}$.*

PROOF. If $d=1$, nothing has to be shown. Otherwise, the first claim is plain if $\mathbf{A}$ is singular. If $\mathbf{A}$ is positive definite, then the claims follow from Lemma A.1(a) with $\mathbf{C} = \mathbf{B} - \mathbf{A}$. □

LEMMA A.3. *Let $\mathbf{x}_1^n$ be a Euclidean data set and let $\{P_1, \ldots, P_g\}$ be a partition of $1..n$ with cardinalities $a_1, \ldots, a_g$. Moreover, let $\mathbf{m}$ be the mean of $\mathbf{x}_1^n$ and let $\mathbf{m}_j$ be the mean of $(\mathbf{x}_i)_{i \in P_j}$ (arbitrary if $a_j = 0$). Then*

$$\sum_{j=1}^{g} \sum_{i \in P_j} (\mathbf{x}_i - \mathbf{m})(\mathbf{x}_i - \mathbf{m})^T$$

$$= \sum_{j=1}^{g} W_{P_j} + \frac{1}{n} \sum_{1 \leq j < h \leq g} a_j a_h (\mathbf{m}_j - \mathbf{m}_h)(\mathbf{m}_j - \mathbf{m}_h)^T.$$

PROOF. Expanding the left-hand side, we obtain

$$\sum_{j=1}^{g} \sum_{i \in P_j} (\mathbf{x}_i - \mathbf{m})(\mathbf{x}_i - \mathbf{m})^T$$

$$= \sum_{i=1}^{n} \mathbf{x}_i \mathbf{x}_i^T - n \cdot \mathbf{m}\mathbf{m}^T$$

$$= \sum_{j=1}^{g} \sum_{i \in P_j} (\mathbf{x}_i \mathbf{x}_i^T - \mathbf{m}_j \mathbf{m}_j^T) + \sum_{j=1}^{g} a_j \mathbf{m}_j \mathbf{m}_j^T - n \cdot \mathbf{m}\mathbf{m}^T$$

$$= \sum_{j=1}^{g} W_{P_j} + \sum_{j=1}^{g} a_j \cdot \mathbf{m}_j \mathbf{m}_j^T - n \cdot \mathbf{m}\mathbf{m}^T.$$

On the other hand,

$$\frac{1}{n} \sum_{1 \leq j < h \leq g} a_j a_h (\mathbf{m}_j - \mathbf{m}_h)(\mathbf{m}_j - \mathbf{m}_h)^T$$

$$= \frac{1}{2n} \sum_{j,h} a_j a_h (\mathbf{m}_j - \mathbf{m}_h)(\mathbf{m}_j - \mathbf{m}_h)^T$$

$$= \frac{1}{n} \sum_{j,h} a_j a_h \mathbf{m}_j \mathbf{m}_j^T - \frac{1}{n} \sum_{j,h} a_j a_h \mathbf{m}_j \mathbf{m}_h^T$$

$$= \sum_{j} a_j \mathbf{m}_j \mathbf{m}_j^T - \frac{1}{n} \left( \sum_{j} a_j \mathbf{m}_j \right) \left( \sum_{j} a_j \mathbf{m}_j \right)^T$$

$$= \sum_{j} a_j \mathbf{m}_j \mathbf{m}_j^T - n \cdot \mathbf{m}\mathbf{m}^T.$$

□

LEMMA A.4. *Let $g \geq 2$ and $c \geq 2$. The minimum of the sum of products $\sum_{1 \leq j < l \leq g} a_j a_l$ taken over all $g$-tuples $(a_1, a_2, \ldots, a_g)$ of real numbers $a_1, a_2 \geq 1$, $a_3, \ldots, a_g \geq 0$ such that $\sum_{j=1}^{g} a_j = c$ is $c - 1$. It is assumed exactly at the tuples $(1, c-1, 0, \ldots, 0)$ and $(c-1, 1, 0, \ldots, 0)$.*

PROOF. We proceed by induction on $g$. If $g = 2$, we have the one-dimensional problem of optimizing the function $a_1 \mapsto a_1(c - a_1)$, restricted to the interval $[1, c-1]$. It plainly attains its minimum at the two endpoints $1$ and $c-1$.

Assume now that the assertion has already been proved up to $g$ and let us prove it for $g+1$. From $\sum_{1 \leq j < l \leq g+1} a_j a_l = a_{g+1}(c - a_{g+1}) + \sum_{1 \leq j < l \leq g} a_i a_j$ and the induction hypothesis, we infer by means of the principle of dynamic optimization,

$$\min_{\substack{\sum_{j=1}^{g+1} a_j = c \\ a_1 \geq 1, a_2 \geq 1}} \sum_{1 \leq j < l \leq g+1} a_j a_l$$

$$= \min_{a_{g+1} \in [0, c-2]} \left( a_{g+1}(c - a_{g+1}) + \min_{\substack{\sum_{j=1}^{g} a_j = c - a_{g+1} \\ a_1 \geq 1, a_2 \geq 1}} \sum_{1 \leq j < l \leq g} a_j a_l \right)$$

$$= \min_{a_{g+1} \in [0, c-2]} (a_{g+1}(c - a_{g+1}) + c - a_{g+1} - 1)$$

$$= \min_{a_{g+1} \in [0, c-2]} (a_{g+1}(c - a_{g+1} - 1) + c - 1).$$

This is a one-dimensional optimization problem for $a_{g+1} \in [0, c-2]$. The minimizer is $0$ and the minimum is $c - 1$. This concludes the proof. □

**Acknowledgments.** We are grateful to the referee for sharing his or her ideas with us. This helped to improve the paper. We also thank the Editor for kindly pointing out the reference Fraley and Raftery (2002) to us.

Fakultät für Mathematik und Informatik
Universität Passau
D-94030 Passau
Germany
e-mail: ritter@fmi.uni-passau.de

TABLE 1
*Fraction of the estimated regular observations whose squared Mahalanobis distances from their relative estimated population centers are greater than $\chi_8^2(\gamma)$, $d=8$. An estimate of the amount of outliers is the fraction $\frac{n-r}{n}$, for which the values shown in the corresponding row match*
*best the theoretical tail probabilities for the $\chi^2$-distribution shown on top*

| | | | $1 - \gamma$ | | | |
|---|---|---|---|---|---|---|
| $\alpha$ | $\beta$ | $\frac{n-r}{n}$ | **0.05** | **0.025** | **0.01** | **0.001** |
| 0.95 | 0.999 | 0 | 0.097 | 0.069 | 0.033 | 0.002 |
| | | 0.05 | 0.066 | 0.045 | 0.015 | 0 |
| | | 0.10 | 0.039 | 0.019 | 0.004 | 0 |
| | | 0.15 | 0.004 | 0 | 0 | 0 |
| 0.99 | 0.999 | 0 | 0.095 | 0.062 | 0.022 | 0.001 |
| | | 0.05 | 0.066 | 0.044 | 0.012 | 0 |
| | | 0.10 | 0.041 | 0.018 | 0.006 | 0 |
| | | 0.15 | 0.003 | 0 | 0 | 0 |
| 0.999 | 0.999 | 0 | 0.094 | 0.055 | 0.016 | 0 |
| | | 0.05 | 0.068 | 0.042 | 0.012 | 0 |
| | | 0.10 | 0.039 | 0.018 | 0.005 | 0 |
| | | 0.15 | 0.003 | 0 | 0 | 0 |
| 0.999999 | 0.999 | 0 | 0.075 | 0.035 | 0.007 | 0 |
| | | 0.05 | 0.068 | 0.024 | 0 | 0 |
| | | 0.10 | 0.045 | 0.023 | 0.004 | 0 |
| | | 0.15 | 0.011 | 0 | 0 | 0 |
| 0.95 | 0.999999 | 0 | 0.100 | 0.094 | 0.086 | 0.051 |
| | | 0.05 | 0.066 | 0.053 | 0.045 | 0.027 |
| | | 0.10 | 0.043 | 0.018 | 0.007 | 0 |
| | | 0.15 | 0.003 | 0 | 0 | 0 |
| 0.99 | 0.999999 | 0 | 0.104 | 0.093 | 0.085 | 0.039 |
| | | 0.05 | 0.069 | 0.055 | 0.044 | 0.031 |
| | | 0.10 | 0.042 | 0.018 | 0.006 | 0 |
| | | 0.15 | 0.003 | 0 | 0 | 0 |
| 0.999 | 0.999999 | 0 | 0.102 | 0.099 | 0.092 | 0.045 |
| | | 0.05 | 0.068 | 0.057 | 0.050 | 0.027 |
| | | 0.10 | 0.041 | 0.019 | 0.005 | 0 |
| | | 0.15 | 0.004 | 0 | 0 | 0 |
| 0.999999 | 0.999999 | 0 | 0.106 | 0.098 | 0.085 | 0.023 |
| | | 0.05 | 0.076 | 0.059 | 0.049 | 0.024 |
| | | 0.10 | 0.046 | 0.023 | 0.004 | 0 |
| | | 0.15 | 0.011 | 0 | 0 | 0 |

TABLE 2
*The maximal Bhattacharyya distances of the best matchings between estimates and theoretical populations*

| | | dimension/number of clusters | | |
|---|---|---|---|---|
| $\alpha$-quantile | $\beta$-quantile | 2/4 | 4/8 | 8/16 |
| 0.95 | 0.99 | — | 0.0685 | 0.0789 |
| | 0.999999 | — | 0.0689 | 0.0556 |
| 0.99 | 0.999 | 0.0386 | 0.0340 | 0.0291 |
| | 0.999999 | 0.0356 | 0.0246 | 0.0297 |
| 0.999 | 0.999 | 0.0257 | 0.0165 | 0.0265 |
| | 0.999999 | 0.0111 | 0.0155 | 0.0265 |
| 0.999999 | 0.999 | 0.0105 | 0.0165 | 0.0241 |
| | 0.999999 | 0.0104 | 0.0176 | 0.0240 |